\documentclass[12pt]{amsart}
\usepackage{amssymb}

\newtheorem{theorem}{Theorem}

\newtheorem{proposition}{Proposition}
\newtheorem{corollary}{Corollary}

\newtheorem{definition}{Definition}
\newtheorem{lemma}{Lemma}
\newcommand{\sm}{\setminus}
\newcommand{\LL}{{\mathcal L}}
\newcommand{\K}{{\mathcal K}}
\newcommand{\KK}{{\mathbb K}}
\newcommand{\KP}{{\mathbb {KP}}}
\newcommand{\Z}{{\mathbb Z}}
\newcommand{\R}{{\mathbb R}}
\newcommand{\PP}{{\mathbb P}}
\newcommand{\RP}{{\mathbb {RP}}}
\newcommand{\C}{{\mathbb C}}
\newcommand{\HH}{{\mathbb H}}
\newcommand{\CP}{{\mathbb {CP}}}

\author{V.A.~Vassiliev}

\email{vva@mi.ras.ru}

\thanks{Supported in part by NWO, project 047--008--005}

\title{Topology of plane arrangements and their complements}

\date{Revised version published in 2001}

\begin{document}

\begin{abstract}
This is a glossary of notions and methods related with the topological theory
of collections of affine planes, including braid groups, configuration spaces,
order complexes, stratified Morse theory, simplicial resolutions, complexes of
graphs, Orlik--Solomon rings, Salvetti complex, matroids, Spanier--Whitehead
duality, twisted homology groups, monodromy theory and multidimensional
hypergeometric functions.

The emphasis on the most geometrical explanation is done; applications and
analogies in the differential topology are outlined.

Some recent results of the theory are presented.
\end{abstract}

\maketitle

\section{Introduction}

Finite collections of affine planes in ${\mathbb R}^N$ or in ${\mathbb C}^N$
(shortly, affine plane arrangements) form a remarkable class of algebraic
varieties. Indeed,

1) they are a meeting point of topology, combinatorics, linear algebra,
representation theory, algebraic geometry, complexity theory, mathematical
physics and differential equations;

2) they are a wonderful proving ground for methods and motivations in these
fields, having very far generalizations;

3) they provide a successful elementary visualization of abstract algebraic and
combinatorial notions and constructions.

Formulas, constructions and theorems once arising in this theory appear then
again and again in very distant fields and problems.

One of main problems of the theory asks to which extent the topological
properties of the union of several planes (and of the complement of this union)
are determined by the formal data, i.e. by the information on the dimensions of
all subcollections of planes.

We shall use this problem to demonstrate such notions and methods as braid
groups, configuration spaces, order complexes, stratified Morse theory,
simplicial resolutions, complexes of graphs, Orlik-Solomon rings, matroids,
Spanier--Whitehead duality, twisted homology groups, monodromy theory,
hypergeometric functions, etc.

There are several very good expositions of the theory of arrangements or
different its aspects, see e.g \cite{GS}, \cite{cart}, \cite{OT}, \cite{Z2},
\cite{varchenko} and especially \cite{bjorner}. An exhaustive survey of
algebraic aspects of the theory of complex hyperplane arrangements is given in
S.~Yuzvinsky's article \cite{Yuz3}.

In this short article, I tried to give an elementary introduction to the
theory, making emphasis on a) the most geometrical aspects and motivations of
the theory, b) the most recent results not reflected yet in introductory texts,
c) the subjects that traditionally are treated in more formal and abstract way
than it is necessary, d) the results having important applications and
generalizations in the fields familiar to me: differential topology,
singularity theory, integral geometry, complexity theory...

\section{Main definitions, notation and examples}
\label{examples}

An {\em affine plane arrangement} is any finite collection of affine planes (of
arbitrary, maybe different, dimensions) in $\R^N$.

An arrangement is called {\em central} if all its planes contain the origin in
$\R^N$. In this case one says also that we have a {\em subspace arrangement}.

One can define also the {\em complex} plane or subspace arrangements in $\C^N$:
they are a special case of usual arrangements in $\R^{2N}$.

In a similar way one defines plane arrangements in $\RP^N$ or in $\CP^N$: they
are in the obvious one-to-one correspondence with central arrangements of
nontrivial subspaces in $\R^{N+1}$ (respectively, in $\C^{N+1}$).

Any real affine plane in $\R^N$ defines a complex plane of the same dimension
in $\C^N$: its {\em complexification}. Therefore the complexification of any
real plane arrangement is well defined.

Suppose that we have a plane arrangement $\LL$ consisting of planes $L_1,
\dots, L_m.$ The union $L_1 \cup \dots \cup L_m$ of these planes is called {\em
the support} of $\LL$ and will be denoted by $L$. For any subset of indices $I
\subset \{1, \dots, m\}$ we set
\begin{equation}
L_I \equiv \bigcap_{i \in I} L_i . \label{inters}
\end{equation}

The first example of a {\em hyper}plane arrangement is provided by the {\em
coordinate cross} in $\R^N$ given by the equation $x_1 \cdot \ldots \cdot x_N
=0.$

The next famous arrangement, the {\em diagonal} arrangement $A(N,2) \subset
\C^N$, consists of $\binom{N}{2}$ hyperplanes $V_{ij} \equiv V_{ji},$ $\{i\ne
j\} \subset \{1, \dots, N\},$ given by equations $x_i =x_j$. The complement of
this arrangement in $\C^N$ can be considered as the {\em $N$th ordered
configuration space} of $\C^1$, i.e the space of all ordered collections of $N$
pairwise different points in $\C^1$.

More generally, for any $k \in [2,N]$ the {\em $k$-equal arrangement} $A(N,k)$
consists of $\binom{N}{k}$ planes $V(i_1, \dots, i_k),$ $1 \le i_1 < \dots <
i_k \le N,$ given by conditions $x_{i_1} = \dots = x_{i_k}.$ We can define such
arrangements both in $\R^N$ and in $\C^N.$

Another generalization of the arrangement $A(N,2)$ is as follows (see
\cite{Brieskorn}). Consider any finite group $W$ of isometries of the Euclidean
space $\R^N$ generated by reflections in several hyperplanes ({\em mirrors}).
(Such groups are well-known, see \cite{Bourbaki}: {\em irreducible} groups of
this type form four infinite series $A_m$ ($m \ge 1$), $C_m$ ($m \ge 2$), $D_m$
($m \ge 4$), and $I_2(p)$, and seven exceptional cases $G_2, F_4, H_3, H_4,
E_6, E_7$ and $E_8$.) Almost all orbits of the action of $W$ in $\R^N$ have one
and the same cardinality. The union of {\em irregular} orbits of smaller
cardinality consists of the mirrors generating the group and their images under
its action. All components of the complement of this union are simplicial cones
({\em Weyl chambers}). The action of $W$ in the space $\R^N$ extends in the
obvious way to an action in its complexification $\C^N$. The union $D_W$ of
irregular orbits of the latter action consists of complexifications of mirrors
and their orbits; it is called the {\em diagonal} of the group $W$.

For instance let the mirrors be all the hyperplanes given by equations
$x_i=x_j,$ $i \ne j.$ Then the group $W$ is isomorphic to the permutation group
$S(N)$, and the corresponding diagonal arrangement coincides with $A(N,2).$
This is the case $A_{N-1}$ of the classification of reflection groups.

A hyperplane arrangement in $\R^N$ or $\C^N$ or $\RP^N$ or $\CP^N$ {\em has
normal crossings} if for any subset $I \subset \{1, \dots, m\}$ the plane $L_I$
either is empty or its {\em codimension} $N -\dim L_I$ is equal to the
cardinality of $I.$

A hyperplane arrangement in $\R^N$ or $\C^N$ is {\em generic} if, being
augmented by the {\em infinitely distant} hyperplane it becomes an arrangement
with normal crossings in $\RP^N$ (respectively, in $\CP^N$).

It is easy to see that for all generic arrangements of $m$ hyperplanes in
$\C^N$ the corresponding triples $(\CP^N, \C^N, L)$ are homeomorphic. Generic
arrangements form an open dense subset in the space of all ordered collections
of $m$ hyperplanes.

\section{Basic example: cohomology rings of pure braid groups}
\label{basbas}

Denote the open manifold $\C^N \setminus A(N,k)$ by $M(N,k)$.

\begin{proposition}[see \cite{FaN}, \cite{ararr}]
For any $N$, the manifold $M(N,2)$ is a $K(\pi,1)$-space, i.e. all its higher
homotopy groups $\pi_2, \pi_3, \dots$ are trivial. \label{fan}
\end{proposition}

Indeed, forgetting the last point of a $N$-point configuration we obtain a
fiber bundle
\begin{equation}
M(N,2) \to M(N-1,2) \label{forget}
\end{equation}
with fiber equal to $\C^1 \setminus \{N-1$ points$\}$. Proposition \ref{fan}
follows by induction from the exact homotopy sequence of this bundle.

The group $\pi_1(M(N,2))$ is called the {\em pure braid group} with $N$
strings. The fundamental group of the similar space of non-ordered sets is
called simply the {\em braid group with $N$ strings}, see \cite{Artin}. For
algebraic and homological properties of these groups see, in particular,
\cite{Lin}, \cite{Arnold-4}, \cite{Gor}, \cite{Sal2}.

A similar statement holds for any finite reflection group.

\begin{theorem}[see \cite{Brieskorn}, \cite{Dimm}]
For any finite group $W$ acting by reflections in $\R^N$, the corresponding
space $\C^N \sm D_W$ (i.e. the union of regular orbits of the complexified
action in $\C^N$) is a $K(\pi,1)$-space. \label{brdel}
\end{theorem}

This theorem was proved by Brieskorn for reflection groups $C_m,$ $D_m,$ $G_2,$
$F_4$, $I_2(p)$; Deligne has proved it in the general case. The groups
$\pi_1(\C^N \sm D_W)$ for these arrangements are called the {\em Brieskorn
braid groups}, see \cite{Brieskorn}.

Theorem \ref{brdel} implies that the cohomology rings of these groups are equal
to these of spaces $\C^N \setminus D_W$. These rings were calculated in
\cite{ararr} for the case $A_{N-1}$ (i.e. that of the arrangement $D_W \equiv
A(N,2)$) and in \cite{Brieskorn} for all other reflection groups.

Moreover, the complement $M_{\R}(N,3)$ of the {\em real} $3$-equal arrangement
also is a $K(\pi,1)$-space for any $N$: this fact was conjectured by
A.~Bj\"orner and proved by M.~Khovanov \cite{hov2}.
\medskip

The calculation of the cohomology ring in the Arnold's case $M(N,2)$ is based
on the same fiber bundle (\ref{forget}).

\begin{proposition}[see \cite{ararr}]
\label{ararrr} The group $H^*(M(N,2))$ is torsion-free and is isomorphic to the
tensor product $H^*(\phi_{N-1}) \otimes H^*(\phi_{N-2}) \otimes \cdots \otimes
H^*(\phi_1)$ where $\phi_i$ is the wedge of $i$ circles. In particular the
Poincar\'e polynomial of $H^*(M(N,2))$ is equal to $(1+t)(1+2t) \cdot \ldots
\cdot (1+(N-1)t)$.
\end{proposition}

Indeed, it is easy to see that the bundle (\ref{forget}) is {\em homologically
simple}, i.e. the fundamental group of the base acts trivially in the homology
of its fiber (which is homotopy equivalent to $\phi_{N-1}$). The spectral
sequence of this bundle stabilizes in the term $E_2$, therefore we have a ring
isomorphism $H^*(M(N,2)) \simeq H^*(M(N-1,2)) \otimes H^*(\phi_{N-1})$. \quad
$\square$
\medskip

The ring structure in this cohomology, also calculated in \cite{ararr}, is as
follows. For any plane $V_{jk}$ of our arrangement denote by $\omega_{jk}$ the
logarithmic form with singularity at this plane, $\omega_{jk} \equiv
\frac{1}{2\pi i} \frac{dz_j -dz_k}{z_j-z_k}$; the integral of it along a closed
path in $M(N,2)$ equals the linking number of this path with the plane
$V_{jk}$.

\begin{proposition}
\label{braimain} For any three different indices $i,j,k \in [1,N]$ the equality
\begin{equation}
\omega_{ij} \wedge \omega_{jk} + \omega_{jk} \wedge \omega_{ki} + \omega_{ki}
\wedge \omega_{ij} = 0 \label{relation}
\end{equation}
holds identically in $M(N,2)$. In particular the similar identity on the
corresponding cohomology classes holds in the ring $H^*(M(N,2)):$ if
$\Omega_{jk}$ is the cohomology class of the form $\omega_{jk}$, then
\begin{equation}
\Omega_{ij} \smile \Omega_{jk} + \Omega_{jk} \smile \Omega_{ki} + \Omega_{ki}
\smile \Omega_{ij} = 0. \label{relat2}
\end{equation}
The integer cohomology algebra $H^*(M(N,2))$ is canonically isomorphic to the
quotient algebra of the exterior algebra formally generated by $\binom{N}{2}$
elements $\Omega_{jk}$ through the ideal multiplicatively generated by left
parts of all possible expressions (\ref{relat2}) with arbitrary $i,j$ and $k$.
\end{proposition}

\section{Orlik--Solomon ring and cohomology of complements of
complex hyperplane arrangements} \label{osos}

A general statement very similar to Proposition \ref{braimain} holds for an
arbitrary complex central hyperplane arrangement.

Let $L = \{L_1, \dots, L_m\}$ be such an arrangement in $\C^N$, whose planes
$L_i$ are given by linear equations $f_i=0.$ A collection of indices $I \subset
\{1, \dots, m\}$ is {\em dependent} if the codimension of $L_I$ is smaller than
the expected value $|I|$ (i.e. the corresponding equations $f_i$ are linearly
dependent). For any dependent set $I = \{i_1 < \dots < i_k\}$ denote by
$\rho(I)$ the rational differential $(k-1)$-form
\begin{equation}
\sum_{j=1}^k (-1)^j\frac{df_{i_1}}{f_{i_1}} \wedge \dots \wedge
\widehat{\frac{df_{i_j}}{f_{i_j}}} \wedge \dots \wedge
\frac{df_{i_k}}{f_{i_k}}.
\end{equation}
It is easy to see that this form is equal to zero in $\C^N \setminus L$.

\begin{theorem}[see \cite{OS}]
\label{os} For any central complex hyperplane arrangement $\LL$ in $\C^N$, the
integral cohomology ring $H^*(\C^N \setminus L)$ is canonically isomorphic to
the quotient algebra of the exterior algebra on generators $\alpha_j$
corresponding to hyperplanes of $\LL$ through the ideal generated by all
elements
\begin{equation}
\label{rel3} \sum_{j=1}^k (-1)^j \alpha_{i_1} \smile \dots \smile
\widehat{\alpha_{i_j}} \smile \dots \smile \alpha_{i_k}
\end{equation}
corresponding to all dependent collections $I=(i_1, \dots, i_k)$. Moreover, the
same ideal is generated by expressions (\ref{rel3}) over only {\em minimal}
dependent subsets $I$.
\end{theorem}

\noindent {\bf Example}. For the arrangement $A(N,2)$ the minimal dependent
sets are exactly the collections of planes of the form $V_{i_1i_2},$
$V_{i_2i_3},$ \dots $V_{i_{q-1}i_q},$ $V_{i_qi_1}$, where $\{i_1, \dots, i_q\}
\subset \{1, \dots, N\}$ is any set of indices, $q \ge 3$. It is easy to see
that any form (\ref{rel3}) defined by such a subset belongs to the ideal
generated by similar forms with $q=3$.
\medskip

A main step towards Theorem \ref{os} was done in \cite{Brieskorn}.

\begin{corollary}[see \cite{Brieskorn}]
The group $H^*(\C^N \setminus L)$ is torsion-free.
\end{corollary}

The case of not central hyperplane arrangements in $\C^N$ can be easily reduced
to that of central hyperplane arrangements in $\C^{N+1}$.

\section{How much the topology of the complement is defined
by the dimensional data: a summary}

Let $\LL$ be an arbitrary affine plane arrangement in $\R^N$, consisting of $m$
planes. Suppose that for any $I \subset \{1, \dots, m\}$ we know the dimension
of the plane $L_I$ (and whether this plane is empty or not). What can be then
said about the topology of $\R^N \setminus L$? Given a topological invariant of
$\R^N \sm L$, is it determined uniquely by these data ?
\medskip

For homology and cohomology groups of $\R^N \setminus L$ the answer to the last
question is positive, see \cite{GM} and \S \ref{gmf} below.

For the stable homotopy type the answer also is positive, see \cite{congr},
\cite{ZZ} and \S \ref{sthoty} below.

For the multiplicative structure in cohomology the answers are as follows.

\noindent {\bf A}. In the most general situation not, see \cite{Z}, \cite{Z2}.

\noindent {\bf B}. For complex arrangements of arbitrary dimensions: yes. For
hyperplane arrangements it follows from the above Orlik-Solomon theorem. For an
arbitrary complex arrangement this was proved in \cite{CP} for rational
cohomology. Then S.~Yuzvinsky \cite{Yuz} proposed an explicit formula for this
rational cohomology multiplication, and finally it was proved \cite{DGM},
\cite{dLS}, that the same formula expresses the multiplication in the integral
cohomology ring, see \S \ref{mumu} below.

\noindent {\bf C}. There is a more general class of real arrangements for which
a large part of the multiplicative structure is determined by the dimensional
data.
\medskip

\noindent {\bf Definition.} The arrangement $\LL$ in $\R^N$ is called a $\ge
2$-{\em arrangement} if for any two its incident planes $L_I \subsetneqq L_J$
we have $\dim J - \dim I \ge 2.$
\medskip

For instance any complex arrangement satisfies this condition.

For such arrangements an explicit formula of the multiplication in the
cohomology was proved in \cite{DGM}, \cite{dLS}. In particular it implies the
following proposition.

Suppose that two planes $L_I, L_J$ of an arrangement are {\em transversal}, i.e
\begin{equation}
\label{transv} N-\dim (L_I \cap L_J) = (N- \dim L_I)+ (N-\dim L_J).
\end{equation}
Let us fix also an orientation of $\R^N$. Then any choice of orientations of
planes $L_I,$ $L_J$ defines in a standard way an orientation of $L_I \cap L_J
\equiv L_{I\cup J}.$ If the orientations of all planes $L_I$ of our arrangement
are fixed, then for any ordered pair of transversal planes $L_I$, $L_J$ we get
a sign $+$ or $-$ indicating whether the fixed orientation of $L_{I \cup J}$
coincides with the orientation defined by fixed orientations of $L_I$ and
$L_J$.

\begin{proposition}
\label{ge2} Suppose that we have two plane $\ge 2$-arrangements $\LL,$ $\LL'$
in $\R^N$ such that $\dim L_I = \dim L'_I$ for any $I$, and there are systems
of orientations of all these planes such that for any pair of multi-indices $I,
J$ satisfying (\ref{transv}) the corresponding signs coincide. Then the
cohomology rings of $\R^N \setminus L,$ $\R^N \setminus L'$ are isomorphic.
\end{proposition}

Certainly, the complex arrangements with equal dimensional data satisfy the
conditions of this proposition: if we choose the complex orientations of all
planes then all indices will be equal to $+$.

\noindent {\bf D}. On the other hand this condition on orientations cannot be
removed: a counterexample see in \cite{Z}.

\noindent {\bf E}. Still, something good can be said even in the most general
case of an arbitrary real affine plane arrangement. The group $H^*(\R^N
\setminus L)$ always admits a natural filtration, see \S \ref{gmf},
\ref{sthoty} below. The corresponding {\em graded ring} is uniquely determined
by the dimensional data and the system of signs as in item C above. Its
description also follows from the results of \cite{DGM}, \cite{dLS}, see \S
\ref{mumu}.
\medskip

The most fragile invariant is the fundamental group of the complement of an
arrangement.

\begin{theorem}[see \cite{R}]
There exist two complex line arrangements in $\C^2$ with equal dimensional data
(i.e. sets of lines having a common point) but with nonisomorhic fundamental
groups of complements of their supports.
\end{theorem}

Finally, if we consider the topology not of the complement of the support $L$
but of this support itself or its one-point compactification $\bar L$, then the
dimensional information is very strong.

\begin{theorem}[see \S \ref{sthoty}]
The homotopy types of topological spaces $L$ and $\bar L$ are completely
defined by the dimensional data.
\end{theorem}

\section{Order complex of a poset. Goresky--MacPherson formula}
\label{gmf}

For an arbitrary real affine plane arrangement, a calculation of the cohomology
group of its complement was given by M.~Goresky and R.~MacPherson as a bright
application of their {\em Stratified Morse theory} \cite{GM}.

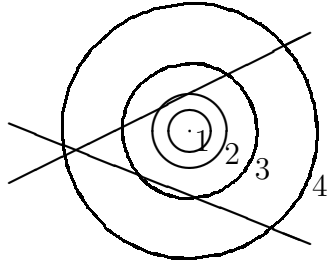
\begin{figure}
\begin{center}
\thicklines \unitlength 1.00mm
\begin{picture}(51.33,36.00)
\put(10.00,20.00){\line(5,-2){40.00}} \put(10.00,12.00){\line(2,1){40.00}}
\put(34.00,19.00){\circle*{0.50}} \put(34.00,19.00){\circle{6.00}}
\put(34.00,19.00){\circle{10.00}}
\multiput(34.00,28.00)(0.98,-0.11){2}{\line(1,0){0.98}}
\multiput(35.97,27.78)(0.31,-0.11){6}{\line(1,0){0.31}}
\multiput(37.84,27.14)(0.19,-0.12){9}{\line(1,0){0.19}}
\multiput(39.53,26.10)(0.12,-0.12){12}{\line(1,0){0.12}}
\multiput(40.95,24.72)(0.11,-0.17){10}{\line(0,-1){0.17}}
\multiput(42.03,23.06)(0.12,-0.31){6}{\line(0,-1){0.31}}
\multiput(42.72,21.21)(0.09,-0.65){3}{\line(0,-1){0.65}}
\multiput(43.00,19.25)(-0.08,-0.99){2}{\line(0,-1){0.99}}
\multiput(42.83,17.27)(-0.12,-0.38){5}{\line(0,-1){0.38}}
\multiput(42.24,15.38)(-0.11,-0.19){9}{\line(0,-1){0.19}}
\multiput(41.25,13.67)(-0.11,-0.12){12}{\line(0,-1){0.12}}
\multiput(39.91,12.21)(-0.16,-0.11){10}{\line(-1,0){0.16}}
\multiput(38.28,11.08)(-0.26,-0.11){7}{\line(-1,0){0.26}}
\multiput(36.45,10.34)(-0.65,-0.11){3}{\line(-1,0){0.65}}
\put(34.50,10.01){\line(-1,0){1.98}}
\multiput(32.52,10.12)(-0.38,0.11){5}{\line(-1,0){0.38}}
\multiput(30.61,10.66)(-0.22,0.12){8}{\line(-1,0){0.22}}
\multiput(28.87,11.60)(-0.14,0.12){11}{\line(-1,0){0.14}}
\multiput(27.38,12.90)(-0.12,0.16){10}{\line(0,1){0.16}}
\multiput(26.21,14.50)(-0.11,0.26){7}{\line(0,1){0.26}}
\multiput(25.41,16.31)(-0.09,0.49){4}{\line(0,1){0.49}}
\put(25.03,18.26){\line(0,1){1.98}}
\multiput(25.09,20.24)(0.10,0.38){5}{\line(0,1){0.38}}
\multiput(25.57,22.16)(0.11,0.22){8}{\line(0,1){0.22}}
\multiput(26.47,23.92)(0.11,0.14){11}{\line(0,1){0.14}}
\multiput(27.72,25.45)(0.14,0.11){11}{\line(1,0){0.14}}
\multiput(29.29,26.67)(0.22,0.11){8}{\line(1,0){0.22}}
\multiput(31.08,27.51)(0.58,0.10){5}{\line(1,0){0.58}}
\multiput(34.00,36.00)(1.09,-0.11){3}{\line(1,0){1.09}}
\multiput(37.27,35.68)(0.39,-0.12){8}{\line(1,0){0.39}}
\multiput(40.41,34.74)(0.22,-0.12){13}{\line(1,0){0.22}}
\multiput(43.32,33.22)(0.14,-0.11){18}{\line(1,0){0.14}}
\multiput(45.87,31.17)(0.12,-0.14){18}{\line(0,-1){0.14}}
\multiput(47.99,28.66)(0.11,-0.20){14}{\line(0,-1){0.20}}
\multiput(49.59,25.79)(0.11,-0.35){9}{\line(0,-1){0.35}}
\multiput(50.60,22.67)(0.10,-0.81){4}{\line(0,-1){0.81}}
\multiput(51.00,19.41)(-0.12,-1.64){2}{\line(0,-1){1.64}}
\multiput(50.76,16.14)(-0.11,-0.40){8}{\line(0,-1){0.40}}
\multiput(49.90,12.97)(-0.11,-0.23){13}{\line(0,-1){0.23}}
\multiput(48.44,10.03)(-0.12,-0.15){17}{\line(0,-1){0.15}}
\multiput(46.45,7.42)(-0.13,-0.11){19}{\line(-1,0){0.13}}
\multiput(43.99,5.25)(-0.20,-0.12){14}{\line(-1,0){0.20}}
\multiput(41.16,3.58)(-0.31,-0.11){10}{\line(-1,0){0.31}}
\multiput(38.07,2.49)(-0.81,-0.12){4}{\line(-1,0){0.81}}
\multiput(34.82,2.02)(-1.64,0.08){2}{\line(-1,0){1.64}}
\multiput(31.54,2.18)(-0.46,0.11){7}{\line(-1,0){0.46}}
\multiput(28.36,2.96)(-0.25,0.12){12}{\line(-1,0){0.25}}
\multiput(25.38,4.35)(-0.16,0.11){17}{\line(-1,0){0.16}}
\multiput(22.73,6.28)(-0.12,0.13){19}{\line(0,1){0.13}}
\multiput(20.49,8.68)(-0.12,0.19){15}{\line(0,1){0.19}}
\multiput(18.76,11.47)(-0.12,0.31){10}{\line(0,1){0.31}}
\multiput(17.60,14.53)(-0.11,0.65){5}{\line(0,1){0.65}}
\put(17.04,17.77){\line(0,1){3.28}}
\multiput(17.12,21.05)(0.12,0.53){6}{\line(0,1){0.53}}
\multiput(17.83,24.25)(0.12,0.27){11}{\line(0,1){0.27}}
\multiput(19.14,27.26)(0.12,0.17){16}{\line(0,1){0.17}}
\multiput(21.01,29.96)(0.12,0.11){20}{\line(1,0){0.12}}
\multiput(23.35,32.25)(0.18,0.12){15}{\line(1,0){0.18}}
\multiput(26.10,34.05)(0.28,0.11){11}{\line(1,0){0.28}}
\multiput(29.14,35.29)(0.81,0.12){6}{\line(1,0){0.81}}
\put(39.67,16.00){\makebox(0,0)[cc]{$2$}}
\put(43.67,14.00){\makebox(0,0)[cc]{$3$}}
\put(51.33,11.67){\makebox(0,0)[cc]{$4$}}
\put(35.67,17.67){\makebox(0,0)[cc]{$1$}}
\end{picture}
\end{center}
\caption{Evolution of sets of smaller values} \label{gm}
\end{figure}

In our special case this theory works as follows. Let us fix a generic point
$X_0 \in \R^N \setminus L$ and a positive quadratic function $f: \R^N \to \R$
with origin at this point (i.e. equal to $x_1^2 + \dots + x_N^2$ in some affine
coordinates centered at $X_0$). The {\em critical} values of $f$ are the
critical values of its restrictions on all nonempty planes $L_I$ of our
arrangement. If $X_0$ and $f$ actually are generic then all these values
corresponding to different planes $L_I$ are different. Further, for any
positive $t$ let us consider the ball $B_t \equiv \{x:f(x) \le t\}$ in $\R^N$
and the {\em manifold of lower values} $\Lambda(t) \equiv B_t \setminus L$. If
the segment $[a,b]$ contains no critical values of $f$ then $\Lambda(a)$ and
$\Lambda(b)$ are homotopy equivalent (and even homeomorphic). If $t$ is
sufficiently small then $\Lambda(t)$ is a ball; if $t$ is sufficiently large
then $\Lambda(t)$ is a deformation retract of the desired space $\R^N \sm L$.
Thus $\R^N \sm L$ can be constructed from a topologically trivial space by a
finite sequence of local surgeries corresponding to all critical values of $f$.

For instance let us consider the line arrangement in $\R^2$ shown in
Fig.~\ref{gm} and its complexification $\LL$ in $\C^2$. There are four
essentially different noncritical values $t_1 < t_2 < t_3 < t_4$: the
intersections of the corresponding balls $B_{t_i}$ with $\R^2$ are shown in
Fig.~\ref{gm}. Let $\Lambda_i \equiv \Lambda(t_i)$ be the corresponding
varieties of our inductive process.

The manifold $\Lambda_1$ is topologically trivial. Passages from $\Lambda_1$ to
$\Lambda_2$ and from $\Lambda_2$ to $\Lambda_3$ are homotopy equivalent to
gluing 1-dimensional cells, so that $\Lambda_3$ is homotopy equivalent to the
bouquet of two circles. Finally, passage from $\Lambda_3$ to $\Lambda_4$ is
equivalent to addition of a 2-dimensional cell. And indeed, it is easy to see
that the resulting manifold $\Lambda_4$ is homotopy equivalent to the
two-dimensional torus.
\medskip

In general, for an arbitrary arrangement the inductive calculation of homology
groups of $\R^N \setminus L$ includes many local problems of the following
sort. Suppose that values $a,b$ are non-critical and the segment $[a,b]$
contains exactly one critical value of $f$; namely, it is the critical value of
the restriction of $f$ to some plane $L_I$. Topological types of manifolds
$\Lambda_a, \Lambda_b$ differ by a surgery localized in a small neighborhood of
the corresponding critical point of $f$ on $L_I$. How does it relate the
cohomology groups of these manifolds?

This problem was solved in \cite{GM} in the combinatorial terms of our
arrangement. The answer is formulated in terms of the {\em order complex} of a
partially ordered set (= {\em poset}).

\begin{definition}
Given a poset $(A,<)$, the corresponding order complex $\Upsilon(A)$ is the
simplicial complex, whose vertices are the points of the set $A$, and the
simplices span all the sequences of such points monotone with respect to the
partial order.
\end{definition}

Every plane arrangement $\LL = \{L_1, \dots, L_m\}$ defines the poset of all
corresponding nonempty sets $L_I,$ $I \subset \{1, \dots, m\}$, and hence the
order complex $\Upsilon(\LL)$.

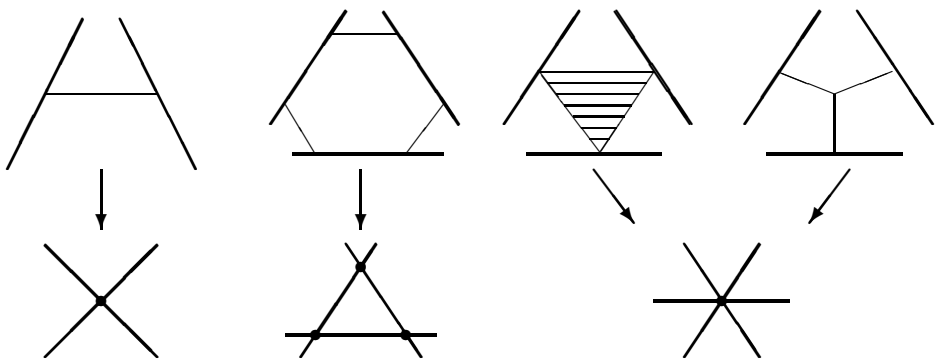
\begin{figure}
\begin{center}
\unitlength 1.00mm \special{em:linewidth 0.4pt}
\thicklines
\begin{picture}(126.67,46.10)
\put(5.10,0.00){\line(1,1){14.90}} \put(5.00,0.10){\line(1,1){14.90}}
\put(5.00,14.90){\line(1,-1){14.90}} \put(5.10,15.00){\line(1,-1){14.90}}
\put(12.50,7.50){\circle*{1.50}} \put(0.10,25.00){\line(1,2){10.00}}
\put(0.00,25.10){\line(1,2){10.00}} \put(25.00,25.10){\line(-1,2){10.00}}
\put(25.10,25.00){\line(-1,2){10.00}} \put(12.50,25.00){\vector(0,-1){8.00}}
\put(39.10,0.00){\line(2,3){10.00}} \put(39.00,0.10){\line(2,3){10.00}}
\put(45.10,15.00){\line(2,-3){10.00}} \put(45.00,15.10){\line(2,-3){10.00}}
\put(57.00,3.05){\line(-1,0){20.00}} \put(57.00,2.95){\line(-1,0){20.00}}
\put(35.10,31.00){\line(2,3){10.00}} \put(35.00,31.10){\line(2,3){10.00}}
\put(50.10,46.00){\line(2,-3){10.00}} \put(50.00,45.90){\line(2,-3){10.00}}
\put(38.00,27.05){\line(1,0){20.00}} \put(38.00,26.95){\line(1,0){20.00}}
\put(53.00,3.00){\circle*{1.50}} \put(41.00,3.00){\circle*{1.50}}
\put(47.00,12.00){\circle*{1.50}} \put(47.00,25.00){\vector(0,-1){8.00}}
\put(90.00,0.10){\line(2,3){10.00}} \put(90.10,0.00){\line(2,3){10.00}}
\put(90.10,15.00){\line(2,-3){10.00}} \put(90.00,15.10){\line(2,-3){10.00}}
\put(86.00,7.45){\line(1,0){18.00}} \put(86.00,7.55){\line(1,0){18.00}}
\put(95.00,7.50){\circle*{1.50}} \put(66.10,31.00){\line(2,3){10.00}}
\put(66.00,31.10){\line(2,3){10.00}} \put(80.77,46.00){\line(2,-3){10.00}}
\put(81.00,46.10){\line(2,-3){10.00}} \put(69.00,27.05){\line(1,0){18.00}}
\put(69.00,26.95){\line(1,0){18.00}} \put(98.10,31.00){\line(2,3){10.00}}
\put(98.00,31.10){\line(2,3){10.00}} \put(113.10,46.00){\line(2,-3){10.00}}
\put(113.00,46.10){\line(2,-3){10.00}} \put(101.00,27.05){\line(1,0){18.00}}
\put(101.00,26.95){\line(1,0){18.00}} \put(112.00,25.00){\vector(-3,-4){5.33}}
\put(78.00,25.00){\vector(3,-4){5.33}} \thinlines
\put(5.00,35.00){\line(1,0){15.00}} \put(53.00,27.00){\line(3,4){5.00}}
\put(41.00,27.00){\line(-3,5){4.00}} \put(43.00,43.00){\line(1,0){9.00}}
\put(79.00,27.00){\line(-3,4){8.33}} \put(70.67,38.00){\line(1,0){15.33}}
\put(86.00,38.00){\line(-2,-3){7.33}} \put(85.00,36.50){\line(-1,0){13.15}}
\put(84.00,35.00){\line(-1,0){11.00}} \put(83.00,33.50){\line(-1,0){8.85}}
\put(82.00,32.00){\line(-1,0){6.75}} \put(81.00,30.50){\line(-1,0){4.65}}
\put(80.00,29.00){\line(-1,0){2.50}} \put(110.00,27.00){\line(0,1){8.00}}
\put(110.00,35.00){\line(-5,2){7.33}} \put(110.00,35.00){\line(5,2){7.67}}
\end{picture}
\end{center}
\caption{Simplicial resolutions of line arrangements} \label{resol}
\end{figure}

\begin{figure}
\unitlength 1.00mm
\begin{center}
\begin{picture}(118.75,39.00)
\put(1.00,4.00){\line(2,3){10.00}} \put(11.00,19.00){\line(2,-3){10.00}}
\put(21.00,4.00){\circle*{1.50}} \put(1.00,4.00){\circle*{1.50}}
\put(11.00,19.00){\circle*{1.50}} \put(11.00,22.00){\makebox(0,0)[cc]{\small
$(12)$}} \put(1.00,1.00){\makebox(0,0)[cc]{\small $(1)$}}
\put(21.00,1.00){\makebox(0,0)[cc]{\small $(2)$}}
\put(48.00,4.00){\line(1,3){5.00}} \put(53.00,19.00){\line(1,-3){5.00}}
\put(58.00,4.00){\line(1,3){5.00}} \put(63.00,19.00){\line(1,-3){5.00}}
\put(68.00,4.00){\line(1,3){5.00}}
\bezier{284}(73.00,19.00)(43.00,39.00)(48.00,4.00)
\put(48.00,4.00){\circle*{1.50}} \put(58.00,4.00){\circle*{1.50}}
\put(68.00,4.00){\circle*{1.50}} \put(73.00,19.00){\circle*{1.50}}
\put(63.00,19.00){\circle*{1.50}} \put(53.00,19.00){\circle*{1.50}}
\put(48.00,1.00){\makebox(0,0)[cc]{\small $(1)$}}
\put(58.00,1.00){\makebox(0,0)[cc]{\small $(2)$}}
\put(68.00,1.00){\makebox(0,0)[cc]{\small $(3)$}}
\put(53.50,22.00){\makebox(0,0)[cc]{\small $(12)$}}
\put(61.50,22.00){\makebox(0,0)[cc]{\small $(23)$}}
\put(75.00,22.00){\makebox(0,0)[cc]{\small $(13)$}}
\put(98.00,4.00){\line(2,3){10.00}} \put(108.00,19.00){\line(0,-1){15.00}}
\put(118.00,4.00){\line(-2,3){10.00}} \put(108.00,19.00){\circle*{1.50}}
\put(108.00,4.00){\circle*{1.50}} \put(118.00,4.00){\circle*{1.50}}
\put(98.00,4.00){\circle*{1.50}} \put(98.00,1.00){\makebox(0,0)[cc]{\small
$(1)$}} \put(108.00,1.00){\makebox(0,0)[cc]{\small $(2)$}}
\put(118.00,1.00){\makebox(0,0)[cc]{\small $(3)$}}
\put(108.00,22.00){\makebox(0,0)[cc]{\small $(123)$}}
\end{picture}
\end{center}
\caption{Order complexes for some line arrangements} \label{ordcom}
\end{figure}
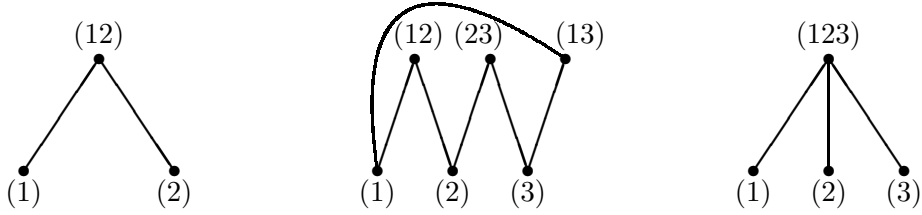

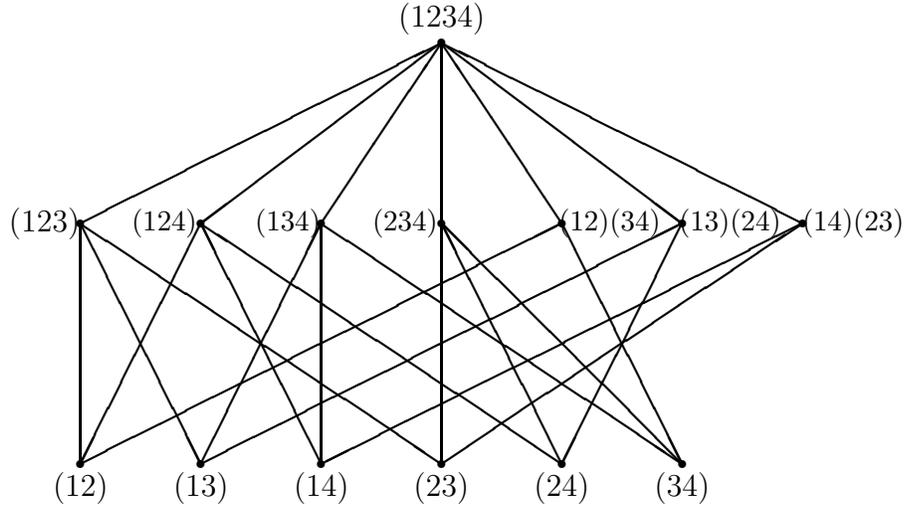
\begin{figure}
\unitlength 0.80mm \special{em:linewidth 0.4pt}
\begin{center}
\begin{picture}(136.50,79.00)
\put(9.00,5.00){\circle*{1.50}} \put(29.00,5.00){\circle*{1.50}}
\put(49.00,5.00){\circle*{1.50}} \put(69.00,5.00){\circle*{1.50}}
\put(89.00,5.00){\circle*{1.50}} \put(109.00,5.00){\circle*{1.50}}
\put(9.00,45.00){\circle*{1.50}} \put(29.00,45.00){\circle*{1.50}}
\put(49.00,45.00){\circle*{1.50}} \put(69.00,45.00){\circle*{1.50}}
\put(89.00,45.00){\circle*{1.50}} \put(109.00,45.00){\circle*{1.50}}
\put(129.00,45.00){\circle*{1.50}} \put(9.00,5.00){\line(0,1){40.00}}
\put(9.00,45.00){\line(1,-2){20.00}} \put(29.00,5.00){\line(1,2){20.00}}
\put(49.00,45.00){\line(0,-1){40.00}} \put(49.00,5.00){\line(-1,2){20.00}}
\put(29.00,45.00){\line(-1,-2){20.00}} \put(9.00,5.00){\line(2,1){80.00}}
\put(89.00,45.00){\line(1,-2){20.00}} \put(109.00,5.00){\line(-1,1){40.00}}
\put(69.00,45.00){\line(0,-1){40.00}} \put(69.00,5.00){\line(3,2){60.00}}
\put(129.00,45.00){\line(-2,-1){80.00}} \put(29.00,5.00){\line(2,1){80.00}}
\put(109.00,45.00){\line(-1,-2){20.00}} \put(49.00,45.00){\line(3,-2){60.00}}
\put(89.00,5.00){\line(-3,2){60.00}} \put(9.00,45.00){\line(3,-2){60.00}}
\put(89.00,5.00){\line(-1,2){20.00}} \put(69.00,75.00){\circle*{1.50}}
\put(23.00,45.00){\makebox(0,0)[cc]{\small (124)}}
\put(43.50,45.00){\makebox(0,0)[cc]{\small (134)}}
\put(63.00,45.00){\makebox(0,0)[cc]{\small (234)}}
\put(97.00,45.00){\makebox(0,0)[cc]{\small (12)(34)}}
\put(117.00,45.00){\makebox(0,0)[cc]{\small (13)(24)}}
\put(137.50,45.00){\makebox(0,0)[cc]{\small (14)(23)}}
\put(3.00,45.00){\makebox(0,0)[cc]{(123)}} \put(9.00,45.00){\line(2,1){60.00}}
\put(69.00,75.00){\line(-4,-3){40.00}} \put(49.00,45.00){\line(2,3){20.00}}
\put(69.00,75.00){\line(0,-1){30.00}} \put(89.00,45.00){\line(-2,3){20.00}}
\put(69.00,75.00){\line(4,-3){40.00}} \put(129.00,45.00){\line(-2,1){60.00}}
\put(9.00,1.00){\makebox(0,0)[cc]{(12)}}
\put(29.00,1.00){\makebox(0,0)[cc]{(13)}}
\put(49.00,1.00){\makebox(0,0)[cc]{(14)}}
\put(69.00,1.00){\makebox(0,0)[cc]{(23)}}
\put(89.00,1.00){\makebox(0,0)[cc]{(24)}}
\put(109.00,1.00){\makebox(0,0)[cc]{(34)}}
\put(69.00,79.00){\makebox(0,0)[cc]{(1234)}}
\end{picture}
\end{center}
\caption{Poset and order complex for the diagonal arrangement $A(4,2)$}
\label{ord4}
\end{figure}

For instance for three line arrangements shown in the lower row of
Fig.~\ref{resol} the corresponding order complexes are given in
Fig.~\ref{ordcom}. For the diagonal arrangement $A(4,2)$ (see \S
\ref{examples}) the order complex is shown in Fig.~\ref{ord4}.

\thicklines

\thicklines

To any $I$ with $L_I \ne \emptyset$ the order subcomplex $\Upsilon(I) \subset
\Upsilon(\LL)$ is associated: this is the union of all simplices in
$\Upsilon(\LL)$ all whose vertices correspond to planes $L_J$ containing $L_I$.
Any such subcomplex $\Upsilon(I)$ is contractible: indeed, all its maximal
simplices have the common vertex corresponding to the plane $\Upsilon(I)$
itself. Denote by $\partial \Upsilon(I)$ the {\em link} of this subcomplex,
i.e. the union of all its simplices not containing its maximal vertex
$\{\Upsilon(I)\}$.

\begin{theorem}[see \cite{GM}]
Suppose that the segment $[a,b]$ contains exactly one critical value of $f$;
let $L_I$ be the corresponding plane. Then

1) $H^i(\Lambda(b),\Lambda(a)) \simeq H_{N-i-\dim L_I -1}(\Upsilon(I), \partial
\Upsilon(I))$;

2) The exact homological sequence of the pair $(\Lambda(b),\Lambda(a))$
splits, i.e. \\
$\tilde H^*(\Lambda(b)) \simeq H^*(\Lambda(b),\Lambda(a)) \oplus \tilde
H^*(\Lambda(a));$ here $\tilde H^*$ denotes the cohomology group reduced modulo
a point.
\end{theorem}

\begin{corollary}[see \cite{GM}]
For an arbitrary affine plane arrangement $\LL$ in $\R^N,$
\begin{equation}
\tilde H^i(\R^N \sm L) \simeq \bigoplus_I H_{N-i-\dim L_I -1}(\Upsilon(I),
\partial \Upsilon(I)), \label{gomac}
\end{equation}
summation over all nonempty planes $L_I$ of the arrangement.
\end{corollary}

In section \ref{sthoty} we shall explain this formula and prove its homotopical
version.
\medskip

\noindent {\bf Example.} Suppose that we have a complex hyperplane arrangement
with normal crossings in $\C^N$. Then the group $H^r(\C^N \sm L)$ is isomorphic
to $\Z^{\lambda(r)}$ where $\lambda(r)$ is the number of sets of indices
$I=\{i_1, \dots, i_r\}$ of cardinality $r$ such that $L_I \ne \emptyset$.
Indeed, for any such set $I$ the corresponding order subcomplex $\Upsilon(I)$
is (the first barycentric subdivision of) an $(r-1)$-dimensional simplex, its
link coincides with the boundary of this simplex, the group $H_*(\Upsilon(I),
\partial \Upsilon(I))$ has unique non-trivial term $\Z$ in dimension $r-1$, and
$\dim_{\R}L_I =2N-2r$. In particular, for the {\em generic} arrangement of $m$
hyperplanes in $\C^N$ we have $H^r(\C^N \sm L) \sim \Z^{\binom{m}{r}}$ for any
$r=0,1, \dots, N$. In the last case a more strong statement holds: the space
$\C^N \setminus L$ is homotopy equivalent to the $N$-skeleton of (the standard
cell decomposition of) the $m$-dimensional torus. Indeed, let us consider the
{\em universal} hyperplane $m$-arrangement $\dag^m$, i.e. the coordinate cross
in $\C^m.$ There exists an affine embedding $\phi: \C^N \to \C^m$ such that the
planes of our generic arrangement in $\C^N$ are preimages of intersections of
$\phi(\C^N)$ with coordinate planes in $\C^m$. The cohomology group of the
space $\C^m \setminus \dag^m \sim T^m$ can be calculated as above by means of
an arbitrary generic real quadratic function $f: \C^m \to \R$. It has exactly
$2^m-1$ critical values corresponding to all planes $L!_I$, and the surgeries
corresponding to the passages through these critical values are homotopy
equivalent to adding the cells of the {\em standard} cell decompositions of the
torus. Now, we can choose our function $f$ with center at a point of the
embedded plane $\phi(\C^N) \subset \C^m$ and in such a way that it grows very
slowly along $\phi(\C^N)$ and very fast in the transversal directions. The
``balls'' $\{x: f(x)\le t \}$ will then look like pancakes spread along
$\phi(\C^N)$. For some $t$, such a ball will intersect all planes $L_I$ with
$|I|\le N$ but do not intersect any planes of smaller dimensions. The
corresponding manifold $\Lambda(t)$ is homotopy equivalent to $\C^N \setminus
L,$ on the other hand it is homotopy equivalent to the $N$-skeleton of the
torus $T^m \sim \C^m \setminus L!$.

The isomorphism (\ref{gomac}), as it follows from its proof in \cite{GM}, is
not canonical: its realization depends on some choices. However it allows one
to define an important increasing filtration in the group $H^*(\R^N \sm L)$:
for any realization of this isomorphism the term $F_i$ of this filtration
corresponds to the sum of terms $H_{*}(\Upsilon(I), \partial \Upsilon(I))$ over
all planes $L_I$ of codimensions $\le i$. This definition already does not
depend on the choices, see \S\S \ref{sthoty}, \ref{cococo}.

\section{Simplicial resolutions and inclusion-exclusion formula.
Mayer--Vietoris spectral sequence and its modifications} \label{sires}

The simplicial resolutions are a far extension of the combinatorial formula of
inclusions and exclusions. They allow us to study effectively the topology of
spaces represented as unions of several subspaces, whose (multiple)
intersections are much easier than their symmetric differences.

First let us demonstrate this method in the simplest discrete situation.
Suppose that a finite set $S$ is represented as the union of finitely many
finite sets $S_i,$ $i=1, \dots, m,$ and we need to find the cardinality of $S$.
To do it we construct the simplicial resolution of $S$. First, we take all sets
$S_i$ separately. If some two sets $S_i, S_j$ have a common point, then we draw
a segment between the corresponding points in the separated copies of $S_i$ and
$S_j$. If some point belongs to the triple intersection $S_i \cap S_j \cap S_k$
then we get three separated points joined by three segments. On the next step,
we add the ``interior part'' of this triangle, then construct tetrahedra over
quadruple intersections, etc. The obtained complex $S'$ is homotopy equivalent
to the initial set $S$: to any point of $S$ there corresponds a simplex in
$S'$. In particular $S$ and $S'$ have equal Euler characteristics. But the
Euler characteristic of the finite set $S$ is its cardinality, while that of
$S'$ is the number of vertices (i.e. of points of all sets $S_i$ taken
separately) minus the number of edges (i.e. the sum of cardinalities of all
sets $S_i \cap S_j$ over all pairs $(i \ne j) \subset \{1, \dots, r\}$) plus
the number of triangles, etc. The result is nothing else than the
exclusion-inclusion formula.

The same method works in the ``continuous'' case, say if the set $S = \cup S_i$
is a $CW$-complex, and all sets $S_i$ and all their intersections $S_I$ are its
cell subcomplexes. Namely, we consider the $(m-1)$-dimensional simplex $\Delta$
whose vertices are in one to one correspondence with the indices $1, \dots, m$.
The simplicial resolution of $S$ can be constructed as a subset in $\Delta
\times S$. For any set of indices $I$ we take the simplex $\Delta(I)$ whose
vertices are the points of the set $I$. The simplicial resolution $S'$ is
defined as the union of all products $\Delta(I) \times S_I$ over all subsets
$I$. The obvious projection $\Delta \times S \to S$ induces the map $S' \to S.$
It is easy to see that this map is proper and is a homotopy equivalence. This
space $S'$ is often much easier to study than the initial space $S.$
\medskip

\noindent {\bf Example.} If there are only two sets $S_1, S_2,$ and $S=S_1 \cup
S_2,$ then we get the Mayer-Vietoris exact sequence
\begin{equation}
\cdots \to H_i(S_1 \sqcup S_2) \to H_i(S) \to H_{i-1}(S_1 \cap S_2) \to
H_{i-1}(S_1 \sqcup S_2) \to \cdots .
\end{equation}

Indeed, the corresponding simplex $\Delta$ is the segment $[1,2].$ The disjoint
union $S_1 \sqcup S_2$ can be realized as the subset $(\{1\} \times S_1) \cup
(\{2\} \times S_2) \subset [1,2] \times (S_1 \cup S_2),$ and we can consider
the exact sequence of the pair $(S', (\{1\} \times S_1) \cup (\{2\} \times
S_2))$.
\medskip

More generally, for an arbitrary number $m$ of sets $S_i$ the resolved complex
$S'$ also has a standard filtration $\phi_1 \subset \cdots \subset \phi_{m}:$
its set $\phi_k$ is the union over $l=1, \dots, k$ of all sets $\Delta(I)
\times S_I$ with $l$-element subsets $I \in \{1, \dots, m\}$. The corresponding
spectral sequence calculating the homology of $S' \sim S$ is called the {\em
Mayer-Vietoris spectral sequence} of the composite set $S = \cup S_i.$ It is
useful in some topological problems, however in the study of plane arrangements
it is quite useless: in this case a different filtration in the resolution set
should be considered, see the next section.

\section{Homotopy type of an affine plane arrangement and
stable homotopy type of its complement} \label{sthoty}

Let our spaces $S_i$ be affine planes in $\R^N$ forming the arrangement $\LL$.
In the bottom row of Fig.~\ref{resol} we give three examples of line
arrangements in $\R^N,$ in the top row their simplicial resolutions are
indicated. Over the right-hand picture we have two different resolutions. The
left one of them is constructed exactly as previously: the preimage of the
central point is the entire simplex $\Delta$ whose vertices correspond to all
lines $L_i$.

In the general situation, let us denote by $\Delta(\LL)$ the union of all faces
$\Delta(I) \subset \Delta$ such that the plane $L_I$ is not empty. The
simplicial resolution of the plane arrangement $\LL$ constructed previously is
a subset of $\Delta(\LL) \times L$. Another, more economical resolution is
constructed as follows. Instead of the complex $\Delta(\LL) \subset \Delta$ we
can consider the order complex $\Upsilon(\LL)$ and define the simplicial
resolution as a subspace of $\Upsilon(\LL) \times L$: namely as the union of
all products $\Upsilon(I) \times L_I$ over all nonempty planes $L_I$ of the
arrangement. The order complexes of three arrangements of Fig.~\ref{resol} are
shown in Fig.~\ref{ordcom}.

For two left arrangements in Fig.~\ref{resol} both constructions give one and
the same space, but for the right-hand one the latter construction gives a
different space, see the very right top picture.

These two constructions are homotopy equivalent. Indeed, the order complex
$\Upsilon(\LL)$ can be identified with a subcomplex of the first barycentric
subdivision of $\Delta(\LL)$: the vertex $\{L_I\}$ of $\Upsilon(\LL)$ goes to
the center of the simplex in $\Delta(\LL)$ spanned by all vertices $\{L_i\}$
with $L_i \supset L_I$. This identification defines an embedding of the latter
construction into the former one. Conversely, for any set of indices $I$ with
nonempty $L_I$ we can send the center of the simplex $\Delta(I)$ into the
vertex $\{L_I\} \in \Upsilon(\LL) \subset \Delta(\LL)$; extending this map by
linearity we obtain a map homotopy inverse to our embedding. We shall call two
constructions of simplicial resolutions using the complexes $\Delta(\LL)$ and
$\Upsilon(\LL)$ the {\em naive} and the {\em economical} simplicial
resolutions, respectively.

Almost all further considerations in this section are equally true for both
constructions. In particular we have the following their properties (see
Proposition \ref{heq}). Suppose that $L_I \ne \emptyset$. Some face $\Delta(J)$
of the simplex $\Delta(I)$ is called {\em marginal} if the corresponding plane
$L_J$ is strictly greater than $L_I$. Denote by $\partial \Delta(I)$ the
subcomplex in $\Delta(I)$ formed by all its marginal faces.

\begin{proposition}
\label{heq} 1. There is a homotopy equivalence $\Delta(\LL) \sim
\Upsilon(\LL).$

\noindent 2. For any nonempty plane $L_I$, the pairs $(\Delta(I), \partial
\Delta(I))$ and $(\Upsilon(I), \partial \Upsilon(I))$ are homotopy equivalent.
\quad $\square$
\end{proposition}

In particular the Goresky--MacPherson formula (\ref{gomac}) can be rewritten in
the following way
\begin{equation}
\label{gomac2} \tilde H^i(\R^N\sm L) \simeq \bigoplus_I H_{N-i-\dim L_I
-1}(\Delta(I), \partial \Delta(I)),
\end{equation}

The space of the resolution of the arrangement with support $L$ will be denoted
by $L'$, and its one-point compactification by $\bar L'$. The next well-known
fact, basic for the entire theory of simplicial resolutions, follows easily
from the Borsuk's lemma, see e.g. \cite{fasis} (unfortunately I do not know the
first reference).

\begin{proposition}
\label{hheq} The obvious projection $\pi: L' \to L$ is a proper map and a
homotopy equivalence. Its extension to a map of one-point compactifications,
$\bar \pi: \bar L' \to \bar L$, also is a homotopy equivalence.
\end{proposition}

\begin{theorem}
\label{vzz} For any finite affine plane arrangement with support $L$, there are
homotopy equivalences
\begin{equation}
\label{abs} L' \sim L \sim \Upsilon(\LL) \sim \Delta(\LL),
\end{equation}
\begin{equation}
\label{rel} \bar L' \sim \bar L \sim \bigvee_I \Sigma^{\dim L_I} (\Upsilon(I)/
\partial \Upsilon(I)) \sim \bigvee_I \Sigma^{\dim L_I} (\Delta(I)/ \partial
\Delta(I)),
\end{equation}
where $\Sigma^k$ denotes the $k$-fold suspension, and the bouquets in
(\ref{rel}) taken over all nonempty planes $L_I$ of the arrangement.
\end{theorem}

The middle equivalence in (\ref{abs}) was found in \cite{GM}, see also
\cite{BLY}.

The middle equivalence $\bar L \sim \bigvee_I \Sigma^{\dim L_I} (\Upsilon(I)/
\partial \Upsilon(I))$ of (\ref{rel}) was proved in \cite{ZZ}; simultaneously
the
 composite equivalence
$$\bar L' \sim \vee_I \Sigma^{\dim L_I}
(\Delta(I)/ \partial \Delta(I))$$ of (\ref{rel}) was proved, see \cite{congr}.
By Propositions \ref{heq}, \ref{hheq} these two equalities involved in the
formula (\ref{rel}) are equivalent to one another and to entire this formula.

For other statements of this type see \cite{nak}, \cite{JOS}.

\begin{corollary}
The Goresky--MacPherson formula $($\ref{gomac}$)$, $($\ref{gomac2}$)$.
\end{corollary}

Indeed, this formula follows from (\ref{rel}) and the Alexander duality
\begin{equation}
\tilde H^i(\R^N \setminus L) \simeq \bar H_{N-i-1}(L)
\end{equation}
where $\bar H_*$ is the {\em Borel--Moore homology group}, i.e. the homology
group of the one-point compactification reduced modulo the added point; cf.
\cite{Arnold-4}. (An equivalent definition: the Borel--Moore homology group is
the homology group of the complex of locally finite singular chains in $X$.)

The resolution space $L'$ has a very useful filtration: its term $F_i(L)$ is
defined as the union of all products $\Delta(I) \times L_I$ (respectively,
$\Upsilon(I) \times L_I$) over all nonempty planes $L_I$ of codimensions $\le
i$. In particular $L' = F_{N} (L).$ This filtration extends to that on $\bar
L'$: the term $\bar F_0$ of the latter filtration is the added point and $\bar
F_i, i \ge 1,$ is the closure of $F_i.$ The filtration mentioned in the end of
\S \ref{gmf} is Alexander dual to the corresponding filtration in the homology
of $\bar L' \sim \bar L$.

\begin{corollary}
\label{co2} The stable homotopy type of the complement of an arbitrary affine
plane arrangement $\LL$ is determined by the dimensions of all its planes
$L_I$, in particular the same is true for all extraordinary homology and
cohomology groups.
\end{corollary}

This corollary is based on the following notion.
\medskip

\noindent {\bf Definition.} Two topological spaces (having homotopy types of
$CW$-complexes) are {\em Spanier--Whitehead dual} to one another if they are
homotopy equivalent to two complementary subsets $X$ and $Y \equiv S^N
\setminus X$ of a sphere $S^N$.
\medskip

The homology and cohomology groups of such spaces are related by the Alexander
duality.

For instance our spaces $\bar L$ and $\R^N \setminus L$ are Spanier--Whitehead
dual to one another.

The important fact (see e.g. \cite{Whit}) is that the Spanier--Whitehead
duality determines an involution on the set of stable homotopy types: all
spaces Spanier--Whitehead dual to stably homotopy equivalent spaces are stably
homotopy equivalent to one another. This reduces Corollary \ref{co2} to Theorem
\ref{vzz}.
\medskip

The use of simplicial resolutions makes the proof of Theorem \ref{vzz}
especially transparent.

Indeed, let us consider the projection $L' \to \Upsilon(\LL)$ induced by the
standard projection of the space $\Upsilon(\LL) \times L \supset L'$. All
fibers of this map are planes $L_I$ for certain sets $I$. The homotopy
equivalence $L' \sim \Upsilon(\LL)$ follows by induction over the consequent
contractions of these fibers over different strata of $\Upsilon(\LL)$, cf.
Lemma 1 in \S III.3.4 of \cite{fasis}; at any step of induction the homotopy
equivalence follows from the Borsuk's lemma.

To prove (\ref{rel}) we use a version of the induction from \S \ref{gmf}. For
any $t>0$ we denote by $\bar L(t)$ the quotient space $L/ (L \cap \{x: f(x) \ge
t\})$. Then for sufficiently small $t$ we have $\bar L(t) =$ \{one point\}; for
sufficiently large $t$ $\bar L(t)$ is homotopy equivalent to $\bar L$. If the
segment $[a,b]$ contains no critical values then $\bar L(a)$ is homotopy
equivalent to $\bar L(b)$, so that all we need is the following lemma.

\begin{lemma}
\label{lemvgm} If the segment $[a,b]$ contains only one critical value of $f$,
namely the critical value of its restriction to the plane $L_I$, then we have a
homotopy equivalence
\begin{equation}
\bar L(b) \sim \bar L(a) \vee \Sigma^{\dim L_I}(\Upsilon(I)/\partial
\Upsilon(I)).
\end{equation}
\end{lemma}

{\em Proof.} Let us consider also the spaces $\bar L'(t) = L'/(L' \cap
\pi^{-1}(\{x:f(x) \ge t\})).$ Then the projection $\pi$ induces homotopy
equivalences $\bar L'(t) \sim \bar L(t),$ and it is enough to prove a version
of Lemma \ref{lemvgm} with $\bar L(b)$ and $\bar L(a)$ replaced by $\bar L'(b)$
and $\bar L'(a)$ respectively. In this proof we use the topological operation
of attaching topological spaces by maps. Namely, given two topological spaces
$X, Y,$ a subspace $A \subset X$ and continuous map $\phi: A \to Y$, the space
$Y \cup_\phi X$ is defined as the quotient space of the disjoint union $X
\sqcup Y$ through the relations $a \sim \phi(a)$ for all $a \in A$. In
particular if we have $Z = X \cup Y$ then $Z$ can be considered as $X$ attached
to $Y$ via the identical embedding $X \cap Y \to Y$. An important property of
this operation is its {\em homotopy invariance}: any homotopy equivalence $f:Y
\to Y'$ induces a homotopy equivalence
\begin{equation}
\label{hominv} Y \cup_\phi X \sim Y' \cup_{f \circ \phi} X.
\end{equation}

For any plane $L_J$ of our arrangement, let $L'_J$ be its {\em proper
preimage}, i.e. the set $\Upsilon(J) \times L_J \subset \Upsilon(\LL) \times
L.$ Let $L!a \subset L'$ be the union of proper preimages of all planes $L_J$
such that $L_J \cap \{x:f(x)<a\} \ne \emptyset$. In particular $L!b = L'_I \cup
L!a.$

Given a subspace $X \subset L',$ denote by $X_{/t}$ its reduction modulo the
set of its points $x$ such that $f \circ \pi(x) \ge t$. Then we have
\begin{equation}
\label{bouk} \bar L'(b) = (L!b)_{/b} = (L'_I)_{/b} \cup (L!a)_{/b} \equiv
(L'_I)_{/b} \cup_{\mbox{id}} (L!a)_{/b}
\end{equation}
where id is the identical embedding $(L'_I)_{/b} \cap (L!a)_{/b} \to \cup
(L!a)_{/b}.$

\begin{lemma}
There is a homotopy equivalence $(L!a)_{/b} \to (L!a)_{/a} \equiv \bar L'(a)$
induced by the reduction modulo the layer $\{x: f \circ \pi(x) \in [a,b]\}$.
\quad $\square$
\end{lemma}

This homotopy equivalence maps the entire set $(L'_I)_{/b} \cap (L!a)_{/b} $
into one point (obtained by the factorization from this layer). Therefore by
(\ref{hominv}) the space (\ref{bouk}) is homotopy equivalent to the wedge of
$\bar L'(a)$ and the quotient space $(L'_I)_{/b}/((L'_I)_{/b} \cap (L!a)_{/b}).
$

The latter space $(L'_I)_{/b}/((L'_I)_{/b} \cap (L!a)_{/b})$ is homotopy
equivalent to $(\Upsilon(I) \times L_I) / ((\Upsilon(I) \times (L_I))_{/b} \cup
(\partial \Upsilon(I) \times L_I)) \sim \Sigma^{\dim L_I} (\Upsilon(I)/
\partial \Upsilon(I))$ (cf. \cite{congr}); Lemma \ref{lemvgm} and Theorem
\ref{vzz} are proved.

\section{Examples: resolutions of important arrangements.
Complexes of connected graphs and hypergraphs}

Let us consider again the diagonal arrangement $A(N,2)$ in $\C^N$ or $\R^N,$
see \S \ref{examples}, and its naive simplicial resolution $L' \subset
\Delta(\LL) \times L$. The smallest plane $L_I$ of this arrangement is the line
$\{x_1 = \cdots = x_N\}$. The preimage of any its point is the entire simplex
$\Delta\equiv \Delta(\LL)$ whose $\binom{N}{2}$ vertices correspond to all
possible hyperplanes $V_{ij} \equiv \{x_i = x_j\}.$

Let us draw somewhere $N$ points labelled by numbers $1, \dots, N.$ It is
convenient to depict any hyperplane $V_{ij}$ by the segment connecting the
points $i$ and $j$. Any face of the simplex $\Delta$ defines the graph
consisting of segments corresponding to all vertices of this simplex. It is
easy to see that the subcomplex of marginal faces $\partial \Delta(\LL) \subset
\Delta(\LL)$ consists of all faces corresponding to not connected graphs. Thus
the homological study of our arrangement appeals to the homology group of the
{\em complex of connected graphs} which is defined as the quotient complex of
the standard (acyclic) triangulation of the simplex $\Delta$ through the
subcomplex spanned by all faces corresponding to all not connected graphs.

\begin{proposition}[see \cite{congr}]
The complex of connected graphs with $N$ vertices is acyclic in all dimensions
other than $N-2$. Its $(N-2)$-dimensional homology group is isomorphic to
$\Z^{(N-1)!}$ and is freely generated by the classes of all {\em snake-like}
(i.e. homeomorphic to a segment) trees, one of whose endpoints is fixed.
\end{proposition}

The first assertion of this proposition is a special case of the Folkman's
theorem on the homology of geometric lattices \cite{Folkman}.
\medskip

\noindent {\bf Remark.} The number $(N-1)!$ already appeared in this work.
Indeed, by Proposition \ref{ararrr} the group $H^{N-1}(M(N,2))$ is
$(N-1)!$-dimensional. In the Goresky--MacPherson formula this group corresponds
to the summand $H_{N-1}(\Upsilon(I),\partial \Upsilon(I)) \equiv$
$H_{N-1}(\Delta(I),\partial \Delta(I))$ where $I$ is the entire set $\{1,
\ldots, N\}$.
\medskip

A nice description of the {\em co}homology group of the same complex is given
in \cite{turdis}.

A natural generalization of this complex is provided by {\em complexes of
$i$-connected graphs}, see \cite{bjo}, \cite{BBLSW}, \cite{tur1}, \cite{tur2}.
They also have important applications in the differential and homotopy
topology.

In a similar way, if we consider the $k$-equal arrangement $A(N,k)$ in $\R^N$
or in $\C^N$ (see \S \ref{examples}) then the smallest plane $L_I$ is again the
line $\{x_1 = \cdots = x_N\}$. The corresponding simplex $\Delta$ has
$\binom{N}{k}$ vertices, its faces are the {\em $k$-hypergraphs} with the same
$N$ nodes $1, \dots, N,$ and the marginal subcomplex $\partial \Delta \subset
\Delta$ is the complex of {\em not connected hypergraphs}. The homology groups
of this complex (and of the complement $M(N,k)$ of the arrangement) were
studied by A.~Bj\"orner and V.~Welker \cite{BW}, see also \cite{Bjorner-2},
\cite{BLY}.

In particular, the following facts were proved.

\begin{theorem}[see \cite{BW}]
For any $k \ge 2$ the simplex with $\binom{N}{k}$ vertices reduced modulo the
union of faces corresponding to non-connected $k$-hypergraphs is homotopy
equivalent to a wedge of spheres, in particular all its homology groups are
torsion-free. Moreover, these groups can be nontrivial only in dimensions equal
to $N-(k-2)t-2$, $1 \le t \le N/k$. The ranks of these groups are multiples of
$\binom{N-1}{k-1}$, and in the higher possible dimension $d=N-k$ the rank of $
H_{N-k}$ is equal to $\binom{N-1}{k-1}$.
\end{theorem}

A general formula for these ranks also is given in \cite{BW} (see Theorem 4.5
there), but it is much more complicated.

\medskip
\noindent {\bf Remark} (see \cite{BLY}). The topology of the real variety
$M(N,k) \subset \R^N$ gives good estimates in the following olympic problem.
Suppose we have $N$ coins, some of which are fake, and a weighbridge. It is
natural to assume that all regular coins are of the same weight, and the
weights of all fake are different from one another and from the weight of
regular coins. Given some $k \ge 2,$ how many measurements is it enough to do
to check that we have at least $k$ regular coins?

Indeed, any measurement separates the space $\R^N$ of all possible collections
of weights into three convex parts: two half-spaces and the hyperplane
separating them. Any sequence of measurements together with their results
specifies some cell (maybe empty) in $\R^N.$ Any correct strategy of solving
our problem (such strategies are called {\em decision trees}) should separate
our space $\R^N$ into such cells, any of which belongs to either the
arrangement $A(N,k)$ or to its complement $M(N,k)$. Thus the number of
generators of the total homology group of either of these spaces provides a
lower estimate of the number of cells, and hence of the complexity of the
decision tree.

Among the origins of the theory of arrangements there is one class of olympic
problems more: that on the cut cake, see \cite{Zas}.

\section{Combinatorial realization of cohomology classes
of complements of arrangements} \label{cococo}

The Goresky--MacPherson formula (\ref{gomac}) has the following direct
realization (found essentially in \cite{ZZ}, the present form given in
\cite{M}). Suppose that an Euclidean metric is fixed in $\R^N$.

Consider a constant vector field $V$ (``power'') in $\R^N$. For any
$r$-dimensional simplex of the order subcomplex $\Upsilon(I)/\partial
\Upsilon(I)$ (i.e. for a strictly decreasing sequence of $r+1$ planes $L_{I_1}
\supset L_{I_2} \supset \ldots \supset L_{I_r} \supset L_I$) and any point $x
\in L_I$ consider the sequence of $r+1$ rays in $\R^N$ issuing from $x$, namely
the trajectories of $x$ in the planes $\R^N, L_{I_1}, \ldots, L_{I_r}$ under
the action of this power. (We can realize $V$ as the gradient field of a linear
function $\theta:\R^N \to \R$, then these rays will be the trajectories of
gradients of restrictions of $\theta$ to these planes.)
\medskip

\noindent {\bf Definition.} The constant vector field $V$ is {\em in general
position} with respect to the plane arrangement $L$ if for any $I$ and any
simplex in $\Upsilon(I)/\partial \Upsilon(I)$ these rays are linearly
independent in $\R^N$.
\medskip

It is easy to see that such vector fields form an open dense subset in the
space $\R^N$ of all constant fields. Let us assume that our field $V$ is
generic. Then for any $I$ and simplex as above the convex hull of our $r+1$
rays is linearly homeomorphic to an $(r+1)$-dimensional octant with origin at
$x$. Such octants over all $x \in L_I$ sweep out an $(r+1+\mbox{dim
}L_I)$-dimensional wedge in $\R^N$.

If we have a $r$-dimensional cycle $\alpha$ of the complex $\Delta(I)/\partial
\Delta(I)$, then the sum of (uniformly oriented) corresponding wedges is a
relative cycle in $\R^N$(mod $L$), and the relative homology class $\nabla
\alpha \in H_{r+1+\dim L_I}(\R^N,L)$ of the latter cycle depends on the class
of $\alpha$ in $H_*(\Delta(I), \partial \Delta(I))$ only (up to a sign
depending on the choice of the orientation of the plane $L_I$).

Finally we take the class in $H^*(\R^N \setminus L)$ Poincar\'e--Lefschetz dual
to $\nabla \alpha $ in $\R^N \setminus L$, i.e. defined by intersection indices
with the relative cycle $\nabla \alpha$.

This realization of the formula (\ref{gomac}) depends on the choice of the
direction $V$, but not very much. Two elements in $\bar H_*(\R^N,L)$,
corresponding in this way to one and the same class $\alpha \in H_*(\Delta(I),
\partial \Delta(I))$ via different generic vector fields can differ by elements
of lower filtration only: more precisely, by sums of similar classes coming
from the summands $H_*(\Delta(J), \partial \Delta(J))$ corresponding to planes
$L_J$ {\it strictly containing} $L_I$.

Moreover, if all planes $L_I$ have codimensions $\ge 2$ in all greater planes
$L_J$, then the isomorphism (\ref{gomac}) is canonical (up to the choice of
orientations of all planes $L_I$): indeed, in this case the space of generic
vectors fields $V$ is path-connected.

By the analogy with the knot theory (cf. \cite{PV}, \cite{vtt}), such
realizations of elements of $H^*(\R^N \setminus L)$ can be called their {\em
combinatorial expressions}.

\medskip
This construction allows one to investigate the multiplicative structure in
$H^*(\R^N \sm L)$, in particular to determine this structure in the associated
graded ring.

\section{Multiplication in cohomology.}
\label{mumu}

Let us rewrite the equality (\ref{gomac}) as that for associated graded groups:

\begin{equation}
\label{gmass} Gr H^*(\R^N \setminus L) \cong \oplus H_{k-*-1-\mbox{\small dim
}L_I}(\Upsilon(I), \partial \Upsilon(I))
\end{equation}

This isomorphism is canonical (up to the choice of orientations of planes
$L_I$), and the multiplication in the associated graded {\it ring} is as
follows.

Let us consider two planes $L_I, L_J \subset L$ and two cycles $A, B$ of the
quotient complexes $\Upsilon(I)/\partial \Upsilon(I)$ and $\Upsilon(J)/\partial
\Upsilon(J)$, dim $A=u$, dim $B=v$, represented by chains (=linear combinations
of simplices) of subcomplexes $\Upsilon(I), \Upsilon(J)$ with boundaries in
$\partial \Upsilon(I)$ and $\partial \Upsilon(J)$ only. The {\em shuffle
product} $A \diamondsuit B$ of these cycles is defined as follows (see
\cite{Yuz}).

If $L_I$ and $L_J$ are not transversal (i.e. belong to some proper plane in
$\R^N$) or have no intersection points, then $A \diamondsuit B =0$. Now suppose
that $L_I$ and $L_J$ are transversal and $L_K= L_I \cap L_J \ne \emptyset$ (we
can take $K=I \cup J$). Let $a \subset A$ and $b \subset B$ be some two
simplices with $u+1$ and $v+1$ vertices respectively, i.e. some decreasing
sequences of intersection planes of $L$ having $\{L_I\}$ and $\{L_J\}$ as their
last elements. Consider all $\binom{u+v+2}{u+1}$ possible {\em shuffles} of
these sequences, i.e. all (non-monotone) sequences of $u+v+2$ planes in which
all elements of $a$ and $b$ appear preserving their orders in the sequences $a$
and $b$. To any such shuffle a monotone sequence corresponds: any element
$\lambda$ of the shuffle coming from the sequence $a$ (respectively, $b$)
should be replaced by the intersection of the corresponding plane with the last
plane coming from the sequence $b$ (respectively, $a$) and staying before
$\lambda$ in the shuffle. The obtained monotone sequence is by definition an
$(u+v+1)$-dimensional simplex of the order complex $\Upsilon(K)$. The shuffle
product of our simplices $a$ and $b$ is defined as the sum of all such
simplices taken with signs equal to parities of the corresponding shuffles
(i.e. numbers of transpositions reducing them to the simple concatenation of
sequences $a$ and $b$) multiplied by one sign more, which depends on
multi-indices $I, J$ and $K$ only and is defined by the comparison of the fixed
coorientation of the plane $L_K$ in $\R^N$ with the ordered pair of
coorientations of $L_I$ and $L_J$. The shuffle product of cycles $A$ and $B$ is
defined by linearity. It is a relative cycle defining an element of the summand
in the right-hand part of (\ref{gmass}) corresponding to the plane $L_K$; this
element depends only on homology classes of $A$ and $B$ in the summands
corresponding to $L_I$ and $L_J$.

\begin{theorem}[see \cite{Yuz}, \cite{Yuz2}, \cite{DGM}, \cite{dLS}]
\label{mult} The isomorphism $($\ref{gmass}$)$ commutes the shuffle product in
its right-hand part and the multiplication in its left part obtained from the
usual cohomological multiplication. If all planes $L_I$ have codimensions $\ge
2$ in all greater planes $L_{I'}$, then the same is true for the isomorphism
$($\ref{gomac}$)$ and the multiplication in the ring $H^*(\R^m \setminus L)$
itself, and not in its graded ring only.
\end{theorem}

This is a corollary of the explicit construction of relative homology classes
described in the previous section. Indeed, the multiplication in the cohomology
ring of an oriented manifold $M$ can be realized as follows. Given two classes
$\alpha, \beta \in H^*(M),$ we take Borel--Moore cycles $[\alpha], [\beta]$
Poincar\'e dual to them and meeting transversally, take their intersection
$[\alpha] \cap [\beta]$ supplied with natural orientation, and consider the
cohomology class Poincar\'e dual to this intersection.

Given two planes $L_I$, $L_J$ of our arrangement and classes $$\alpha \in
H_*(\Upsilon(I), \partial \Upsilon(I)), \qquad \beta \in H_*(\Upsilon(J),
\partial \Upsilon(J)),$$ we can realize corresponding elements in the left part
of (\ref{gomac}) with the help of different constant vector fields $V_I$, $V_J$
in $\R^N$ that are in general position to one another if $L_I$ and $L_J$ have
nonempty transversal intersection; if not then these directions should be
opposite to one another and transversal to a hyperplane separating or
containing these planes.
\medskip

Proposition \ref{ge2} follows immediately from Theorem \ref{mult}.
\medskip

\noindent {\bf Exercise}: deduce the (Orlik--Solomon) Theorem \ref{os} from
this one.

\section{Salvetti complex for complexified real hyperplane arrangement}
\label{salvet}

Consider the complexification of a real affine hyperplane arrangement, i.e. the
set of complex hyperplanes $L_j \subset \C^N,$ $j =1, \dots, m,$ distinguished
by the same linear equations $f_j(x)=0,$ $f_j(\R^N)=\R.$ Its complement $\C^N
\setminus L$ is an $N$-dimensional Stein manifold, in particular is homotopy
equivalent to a cell complex of dimension $\le N$. M.~Salvetti \cite{Sal},
following some ideas of \cite{Dimm}, has constructed explicitly an
$N$-dimensional simplicial complex, embedded into the space $\C^N \sm S$ as its
deformation retract. Here we give an easy description of this construction in
the terms of the dual complex.

For any one of our planes $L_j$, its complement $\C^N \sm L_j$ can be
subdivided into four cells $+_j$ $-_j$, $\uparrow_j$ и $\downarrow_j$, given by
conditions ${\rm Re}f_j>0,$ ${\rm Re}f_j<0,$ $\{{\rm Re}f_j=0, {\rm Im}f_j>0\}$
and $\{{\rm Re}f_j=0, {\rm Im}f_j<0\}$ respectively. To any of $4^m$ possible
sequences of $m$ signs $+,-,\uparrow$ and $\downarrow$ we associate the
intersection of corresponding cells (like e.g. $(+_1) \cap (\uparrow_2) \cap
(-_3) \cap \ldots \cap (\uparrow_m)$). This intersection of several real affine
planes and open half-spaces in $\C^N$ is homeomorphic to a cell. By definition
it lies in $\C^N \setminus L$, and any point of $\C^N \setminus L$ belongs to
exactly one cell of this sort.

\begin{lemma}
\label{7.4.2.lem} The subdivision of the manifold $\C^N \sm L$ into cells
corresponding to all possible sequences of signs $+,-,\uparrow$ and
$\downarrow$, augmented with one 0-dimensional cell, defines a cellular
structure on the one-point compactification of this manifold.
\end{lemma}

The proof is elementary, cf. \cite{Fuchs}. \quad $\square$
\medskip

The {\em Salvetti complex} (as a topological space) is just the complex dual to
this cell decomposition. As a combinatorial object, it is defined as a certain
subdivision of this dual complex.

This construction was used in \cite{GR} for defining some topological
invariants of abstract {\em oriented matroids} (see \S \ref{matrs} and
\cite{ormat}). See also \cite{BZ}.

The above described cell decomposition of $\C^N \sm L$ can be simplified very
much if our arrangement has only normal crossings. In this case to any plane
$L_I$ distinguished by several equations $f_i =0,i \in I,$ $f_i(\R^N) \subset
\R,$ we associate the {\em imaginary wedge} in $\C^N$ distinguished by the
conditions $\{\mbox{Re }f_i =0, \mbox{Im }f_i >0\},$ $i \in I.$ Denote by
$\nabla_I$ this wedge from which its intersections with $L$ and with all
smaller wedges $\nabla_J, $ $J \supsetneqq I,$ are removed.

\begin{lemma}[see \cite{VGZ}]
If $L$ is the complexification of a real hyperplane arrangement with normal
crossings, then any nonempty set $\nabla_I$ is homeomorphic to a cell of
dimension $2N-\# I$, and all these sets form a cellular decomposition of the
one-point compactification of $\C^N \sm L$.
\end{lemma}

All the incidence coefficients of the corresponding cell complex are trivial,
therefore the Borel--Moore homology group $\bar H_*(\C^N \sm L)$ is free
Abelian, with the rank of $\bar H_{2N-p}$ equal to the number of nonempty
planes $L_I$ with $\# I=p.$ Of course, the last statement follows also from the
Goresky--MacPherson formula, but the cones $\nabla_I$ provide an especially
easy its realization.

On combinatorial and topological properties of hyperplane arrangements see also
\cite{var-}, \cite{varchenko}, \cite{VS}, \cite{VGZ}, \cite{Gelfand 86},
\cite{Gelfand-Gelfand 86}, \cite{GZ 86}, \cite{GR}, \cite{GS},
\cite{bjorner}--\cite{BZ}, \cite{GKZ}, \cite{Hattori}, \cite{mnev2},
\cite{nak}, \cite{OS}, \cite{OT}, \cite{Zas}--\cite{ZZ}.

\section{Homology of complements of arrangements with twisted coefficients.
Resonances} \label{twist}

Things become slightly more difficult if we consider homology groups with
coefficients in local systems.
\medskip

\noindent {\bf Definition.} A {\em linear local system} on a (locally simply
connected) topological space $X$ (say, $X=\C^N \sm L$) is a vector bundle
$\pi:M \to X$ with fiber $\C^1$ supplied with a flat connection respecting the
$\C$-module structure in the fibers.
\medskip

In other words, for any point $x \in X$, any sufficiently small neighborhood
$U$ of $x$ in $X$ and any point $a \in \pi^{-1}(x)$ we have a distinguished
section of our bundle over entire $U$ equal to $a$ at $x$; this section will be
the same if we start from some other point $x' \in U$ and the intersection
point $a'$ of the old section with the fiber $\pi^{-1}(x')$. For any two points
$a_1, a_2$ of $\pi^{-1}(x)$ the distinguished section equal to $a_1 + a_2$ at
$x$ consists of fiber-vise sums of images of sections equal to $a_1$ and $a_2$
at $x$, and for any $\lambda \in \C$ the section equal to $\lambda a$ at $x$
consists of multiplied by $\lambda$ images of the section equal to $a$ at $x$.

An $i$-dimensional {\em singular simplex} of the local system $\Theta$ is a
continuous map of a standard simplex $\Delta^i$ into the total space $M$ of the
fiber bundle, respecting the flat connection: if a point $\xi$ of the simplex
goes to some point $a \in \pi^{-1}(x), x \in X,$ then some small neighborhood
of $\xi$ goes into the image of the corresponding section over a small
neighborhood of $x$. The group of $i$-dimensional singular chains with
coefficients in the local system $\Theta$ is defined as the quotient group of
the Abelian group generated by all such locally horizontal maps $\Delta^i \to
M$ through the following conditions:

a) if we have two simplices $\phi_1, \phi_2: \Delta^i \to M$ with the same
projection (i.e. $\pi \circ \phi_1 \equiv \pi \circ \phi_2$) then the sum
$\phi_1 + \phi_2$ of them is equal to the third simplex $\phi_3$ such that
$\phi_3(\xi) = \phi_1(\xi)+\phi_2(\xi)$ for any $\xi \in \Delta^i$;

b) for any $\lambda \in \C$ the singular simplex $\phi$ taken with coefficient
$\lambda$ is equal to the singular simplex mapping any point $\xi$ to $\lambda
\phi(\xi)$.

The boundary of such a singular simplex is defined in the obvious way and is a
sum of singular simplices (of reduced dimension) of the same local system. We
can consider the complex of finite chains (i.e. finite linear combinations of
singular simplices) or locally finite chains (i.e. locally finite combinations
whose projections to $X$ also are locally finite chains there). The
corresponding homology groups will be denoted by $H_i(X,\Theta)$ and
$H_i^{lf}(X,\Theta),$ respectively; they are called homology groups {\em of}
(or {\em with coefficients in}) the local system $\Theta$.

These groups are connected by the canonical homomorphism
\begin{equation}
H_i(X,\Theta) \to H_i^{lf}(X,\Theta). \label{canon}
\end{equation}

In particular if our bundle $M$ is the direct product $X \times \C^1$ with the
obvious flat connection then we get the usual homology groups $H_i(X,\C)$ and
the Borel--Moore homology groups $\bar H_i(X,\C)$ respectively. The cohomology
groups $H^i(X,\Theta),$ $H^i_{lf}(X,\Theta)$ are defined in the standard way as
homology groups of conjugate complexes.

Any linear local system defines in the obvious way the monodromy homomorphism
of the fundamental group $\pi_1(X)$ (or, equivalently, of $H_1(X)$) into the
group $\mbox{Aut} (\C^1) \equiv \C^*$: extending our sections over a closed
loop in $X$ we multiply the fiber by the {\em monodromy coefficient}
corresponding to this loop.

To any local system $\Theta$ there corresponds its {\em dual system}
$\Theta^*$, whose fibers are identified with spaces of $\C$-linear functions on
the fibers of the initial system, and this identification respects the flat
connections in both. The monodromy coefficients defined by one and the same
element of $\pi_1(X)$ in dual local systems are inverse to one another (i.e.
their product is equal to 1). If the space $X$ is a $d$-dimensional oriented
manifold then the {\em Poincar\'e isomorphisms}
\begin{equation}
\label{poinc} H_i^{lf}(X,\Theta) \simeq H^{d-i}(X,\Theta^*), \qquad
H^i_{lf}(X,\Theta) \simeq H_{d-i}(X,\Theta^*)
\end{equation}
relate its homology and cohomology groups with coefficients in dual local
systems. (Moreover, a special local system $\mbox{Or}$, called the {\em
orientation sheaf} makes sense of the Poincar\'e isomorphisms on non-oriented
manifolds: this is true even if $\Theta$ and $\Theta^*$ are the constant local
system. Namely, in any of two equations (\ref{poinc}) the right-hand term
should be replaced by the similar term in which the coefficient local system is
not $\Theta^*$ but $\Theta^* \otimes \mbox{Or}$: its monodromy coefficients
coincide with these of $\Theta^*$ up to a sign $+$ or $-$ depending on whether
the corresponding loops preserve the orientation of $X$ or not.)

Local systems and their homology groups are an adequate tool for the calculus
of ramified differential forms and their integrals, see \cite{Deligne 70}.
Suppose that we have a closed analytic $i$-form on $X$, such that the analytic
continuation along any closed path $c$ in $X$ multiplies it by a complex
number, and this number $\tau(c)$ depends on the class of our path $c$ in the
group $H_1(X)$ only. (An important class of such forms will be considered in \S
\ref{hgf}.) Then the integration cycles for this differential form are well
defined as elements of the $i$-dimensional homology group of $X$ with
coefficients in a local system, whose monodromy coefficient at any path $c$ is
equal to $1/\tau(c)$.
\medskip

Now suppose that we have a hyperplane arrangement $\LL = \{L_1, \dots, L_m\}$
in $\C^N$, $m \ge N,$ $X = \C^N \sm L$, and a linear local system $\Theta$ over
$X$; let us denote by $\tau_1, \dots, \tau_m$ the monodromy coefficients of
this system corresponding to small circles going around these hyperplanes in
the positive direction.

\begin{theorem}[see \cite{VGZ}, \cite{AVGL}, \cite{ramif}]
\label{reson1} Let $\LL$ be a {\em generic} hyperplane arrangement in $\C^N$.
Then

A. If there is at least one coefficient $\tau_i$ not equal to 1, then the
groups $H_i(X,\Theta)$, $H_i^{lf}(X,\Theta)$ are nontrivial only for $i=N$ and
their dimensions are equal to $\binom{m-1}{N}$.

B. The map (\ref{canon}) between these groups is an isomorphism if and only if
all numbers $\tau_i$ $(i=1, \dots, m),$ and their product $\tau_0 \equiv \tau_1
\cdot \ldots \cdot \tau_m$ are not equal to 1.

C. If all these numbers are different from $1$ and the arrangement $\LL$ is the
complexification of a real one, then the group $H_N^{lf}(\C^N \sm L, \Theta)$
is freely generated by classes of all bounded components of $\R^N \sm L$.
\end{theorem}

In the last case of a complexified real generic hyperplane arrangement these
facts (the dimensions of both groups, the bijectivity of the map (\ref{canon}),
and the assertion of item C) were proved first by K.~Aomoto \cite{Aomoto 77} in
much stronger assumptions on the coefficients $\tau_j.$

For non-generic arrangements the set of exceptional values of ${\bf \tau}
\equiv (\tau_1, \ldots, \tau_m)$ for which the map (\ref{canon}) is not
bijective, is more complicated, see in particular Theorem \ref{reson2} below.
Since \cite{VGZ} such values are called {\em resonances} of our local system.
\medskip

\begin{figure}
\begin{center}
\unitlength 1.00mm
\begin{picture}(40.00,21.00)
\put(18.00,1.00){\circle*{1.00}} \put(27.00,5.00){\circle*{1.00}}
\put(5.00,8.00){\circle*{1.00}} \put(37.00,1.00){\circle*{1.00}}
\put(27.00,5.00){\vector(0,1){16.00}} \put(18.00,1.00){\vector(0,1){20.00}}
\put(5.00,8.00){\vector(0,1){13.00}} \put(37.00,1.00){\vector(0,1){20.00}}
\put(8.00,8.00){\makebox(0,0)[cc]{$L_1$}}
\put(8.00,16.00){\makebox(0,0)[cc]{$\nabla_1$}}
\put(41.00,2.00){\makebox(0,0)[cc]{$L_m$}}
\put(41.00,17.00){\makebox(0,0)[cc]{$\nabla_m$}}
\end{picture}
\end{center}
\caption{Decomposition of the complement of a point arrangement in $\C^1$}
\label{rays}
\end{figure}
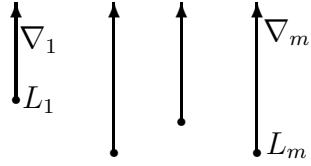

\noindent {\bf Example.} Let be $N=1,$ and $\LL$ the collection of $m$
different points $L_1, \dots, L_m$. If at least one of corresponding
coefficients $\tau_i$ is different from 1 then the group $H^{lf}_1(\C^1 \sm L,
\Theta)$ is generated by the classes of any $m-1$ of $m$ parallel rays
$\nabla_j$ connecting these points with the infinity, see Fig.~\ref{rays}.
Indeed, all these $m$ rays and the complement of the union of them are the
cells covering entire $\C^1 \sm L.$ We can calculate the group $H^{lf}_*(\C^1
\sm L,\Theta)$ with the help of this cellular structure. It is easy to see that
for some natural choice of generators of this cellular complex (i.e. of pairs
\{a cell, a distinguished section of the local system over it\}) the incidence
coefficients of the unique 2-dimensional cell with all 1-dimensional ones are
equal to $\tau_1-1, \dots, \tau_m-1.$ This calculates the group $H^{lf}_*(\C^1
\sm L, \Theta)$; the structure of the ``absolute'' group $H_*(\C^1 \sm
L,\Theta)$ follows by the Poincar\'e duality from the similar statement for the
group $H_*(\C^1 \sm L,\Theta^*)$.

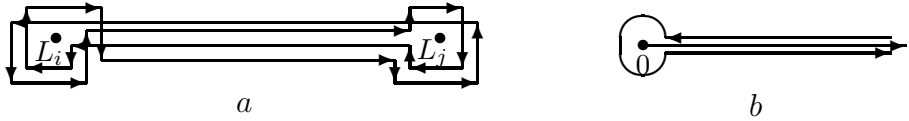
\begin{figure}
\begin{center}
\unitlength 1.00mm
\begin{picture}(120.00,14.00)
\put(7.00,10.00){\circle*{1.50}} \put(58.00,10.00){\circle*{1.50}}
\put(11.00,11.00){\vector(1,0){43.00}} \put(54.00,11.00){\vector(0,1){3.00}}
\put(54.00,14.00){\vector(1,0){7.00}} \put(61.00,14.00){\vector(0,-1){8.00}}
\put(61.00,6.00){\vector(-1,0){7.00}} \put(54.00,6.00){\vector(0,1){3.00}}
\put(54.00,9.00){\vector(-1,0){45.00}} \put(9.00,9.00){\vector(0,-1){3.00}}
\put(9.00,6.00){\vector(-1,0){6.00}} \put(3.00,6.00){\vector(0,1){8.00}}
\put(3.00,14.00){\vector(1,0){10.00}} \put(13.00,14.00){\vector(0,-1){7.00}}
\put(13.00,7.00){\vector(1,0){39.00}} \put(52.00,7.00){\vector(0,-1){3.00}}
\put(52.00,4.00){\vector(1,0){11.00}} \put(63.00,4.00){\vector(0,1){8.00}}
\put(63.00,12.00){\vector(-1,0){62.00}} \put(1.00,12.00){\vector(0,-1){8.00}}
\put(1.00,4.00){\vector(1,0){10.00}} \put(11.00,4.00){\vector(0,1){7.00}}
\put(6.00,8.00){\makebox(0,0)[cc]{{\small $L_i$}}}
\put(57.00,8.00){\makebox(0,0)[cc]{{\small $L_j$}}}
\put(85.00,9.00){\circle*{1.50}} \put(85.00,9.00){\vector(1,0){35.00}}
\put(118.00,10.00){\vector(-1,0){30.00}} \put(85.00,10.00){\oval(6.00,6.00)[t]}
\put(85.00,8.00){\oval(6.00,6.00)[b]} \put(82.00,10.00){\line(0,-1){2.00}}
\put(88.00,8.00){\vector(1,0){30.00}} \put(100.00,1.00){\makebox(0,0)[cc]{$b$}}
\put(32.00,1.00){\makebox(0,0)[cc]{$a$}}
\put(85.00,6.50){\makebox(0,0)[cc]{{\small $0$}}}
\end{picture}
\end{center}
\caption{The cycle ``double loop''} \label{doloop}
\end{figure}

If $\tau_j = 1$ then a small circle around the point $L_j$ is a nontrivial
element of the group $H_1(\C^1 \sm L)$ (because its intersection index with the
ray $\nabla_j$ is not equal to zero). On the other hand it obviously is
homological to zero via a locally finite 2-chain bounded by it, therefore the
kernel of the map (\ref{canon}) is nontrivial. If $\tau_1 \cdot \ldots \cdot
\tau_m =1$ then the same is true for the big circle embracing all points $L_j$.

If all numbers $\tau_1, \dots, \tau_m,$ and their product $\tau_0 \equiv \tau_1
\cdot \ldots \cdot \tau_m$ are different from 1 then the operator inverse to
(\ref{canon}) is provided by {\em double loops}, i.e. cycles shown in
Fig.~\ref{doloop} a). Such a cycle in $\C^1 \sm L$ goes twice (and in opposite
directions) around any of our two singular points $L_i$, $L_j$, therefore it
can be lifted to a cycle of the local system $\Theta$ and defines an element of
$H_1(\C^1 \sm L)$. Since $(1-\tau_i)(1-\tau_j) \ne 0$ the image of this cycle
in $H_1^{lf}(\C^1 \sm L)$ is equal to the class of the interval $(L_i,L_j)$
taken with some nonzero coefficient.
\medskip

Now let us prove Theorem \ref{reson1} in the case of arbitrary $N$.

A. Consider again the universal hyperplane $m$-arrangement, i.e. the coordinate
cross $\dag^m \subset \C^m$, and the local system $\tilde \Theta$ on $\C^m \sm
\dag^m$ with the same monodromy coefficients $\tau_j$. If $\tau_j \ne 1$ for at
least one $j$ then the group $H_*(\C^m \sm \dag^m,\tilde \Theta)$ is trivial in
all dimensions by the K\"unneth theorem in $(\C^m \sm \dag^m) = (\C^*)^m$. Our
generic arrangement $L \subset \C^N$ can be realized as the preimage of
$\dag^m$ under some generic affine embedding $\C^N \to \C^m;$ our local system
$\Theta$ coincides with one induced from $\tilde \Theta$ by the same embedding.
By the Lefschetz theorem (see \cite{GM}), our embedding $\C^N \sm L \to \C^m
\sm \dag^m$ is $N$-connected, in particular induces an isomorphism of homology
groups in all dimensions lower than $N$. Thus $H_i(\C^N \sm L, \Theta)$ is
trivial for $i<N$. By the Poincar\'e isomorphism (\ref{poinc}) the same is true
for all groups $H_i^{lf}(\C^N \sm L, \Theta)$ with $i>N$. The dimension of
these groups in the unique remaining dimension $N$ follows from the
considerations with the Euler characteristic (which does not depend on the
choice of the system $\Theta$).

B. If at least one of numbers $\tau_j$ is equal to $1$ then, similarly to the
one-dimensional example, the small circle around the corresponding plane $L_j$
is a nontrivial element of the kernel of the map (\ref{canon}). If all numbers
$\tau_j$ and their product $\tau_0$ are different from 1 then the bijectivity
of (\ref{canon}) follows from a much more general fact.

\begin{proposition}
\label{gen} Let $W$ be a compact $N$-dimensional complex algebraic manifold,
and $\tilde L$ a finite set of smooth divisors in it having normal crossings
only; let $\theta$ be a linear local system on $W \sm \tilde L$ such that all
the monodromy coefficients corresponding to small circles going around
components of $\tilde L$ are different from 1. Then the canonical map
$$H_*(W \sm \tilde L,\theta) \to
H_*^{lf}(W \sm \tilde L,\theta)$$ is an isomorphism in all dimensions.
\end{proposition}

In particular if $W \sm \tilde L$ is a Stein manifold then both groups can be
nontrivial only in the dimension $N$. This proposition follows easily from the
{\em Leray spectral sequence} for sheaf cohomology, see e.g. \cite{AVGL},
\cite{ramif}.

Our assertion on bijectivity of (\ref{canon}) follows from this proposition if
we take $W = \CP^N$ and $\tilde L=$ the union of $L$ and the improper plane.

C. If our arrangement is the complexification of a real one, then the reversing
of the map (\ref{canon}) can be visualized by the ``multidimensional double
loops'' generalizing Fig.~\ref{doloop} a). The construction of them was
announced in \cite{VGZ} and described in \cite{ramif}. Namely, let $\Delta
\subset \R^N$ be a bounded connected component of $\R^N \setminus L$. The
corresponding ``double loop'' is an $N$-dimensional manifold $\kappa(\Delta)$
together with an immersion $K: \kappa(\Delta) \to \C^N \sm L$ such that

a) this immersion can be lifted to a map into the space $M$ of our local system
$\Theta$, locally flat with respect to its connection; in particular it defines
(up to a scalar coefficient depending on the choice of this lifting) an element
of the group $H_N(\C^N \sm L, \Theta)$;

b) the map (\ref{canon}) sends this element into the homology class of the
component $\Delta$ taken with some coefficient, which is different from zero if
and only if all monodromy coefficients $\tau_j$ corresponding to all walls of
our polytope $\Delta$ are not equal to 1.

This cycle can be constructed with the help of the embedding of $(\C^N,L)$ in
the space $(\C^p,\dag^p)$ of the universal $p$-arrangement, where $p$ is the
number of walls of $\Delta$. First, we construct an immersion $\R^1 \to \C^1
\sm \{0\}$ shown in Fig.~\ref{doloop} b): it coincides with the identical map
in the segment $[+\epsilon, +\infty)$, with the map $\{x \to -x\}$ in the
segment $(-\infty, -\epsilon]$, and maps the segment $[-\epsilon, +\epsilon]$
to a small loop around the origin. The direct product of $m$ copies of such
immersions defines an immersion $\Xi_p: \R^p \to \C^p \sm \dag^p$ that covers
the positive octant in $\R^p \subset \C^p$ with multiplicity $2^p$.

Let $\psi_i, i=1, \dots, p,$ be the linear functions $\C^N \to \C$
distinguishing the planes $L_j$ bounding the component $\Delta$; they can be
normed so that all of them take positive values in $\Delta$. These functions
$(\psi_1, \dots, \psi_p)$ define an embedding $\Psi: \C^N \to \C^p$. The
desired cycle $\kappa(\Delta)$ is induced by this map from the universal
immersion $\Xi_p$. Namely, $\kappa(\Delta)$ is defined as the subset in the
direct product $\R^p \times (\C^N \sm L)$ consisting of such pairs $(t,x)$ that
$\Xi_p(t) = \Psi(x)$. The immersion $K$ is the restriction of the obvious
projection $\R^p \times (\C^N \sm L) \to \C^N \sm L$.
\medskip

\noindent {\bf Remark.} Similar ``double loops'' reverse the map (\ref{canon})
also in the most general situation described in Proposition \ref{gen}, but
their construction is much more complicated, see \S I.10 in \cite{ramif}.
(Direct construction of similar cycles in some interesting particular cases was
given in \cite{Pham 65}.) In Calculus this reversion is called the {\em
regularization of improper integrals}, see also \S \ref{hgf} below.
\medskip

Let us consider a generic linear function $f : \R^N \to \R$ whose values at all
$0$-dimensional planes $L_J$ are different, and associate with any bounded
component $\Delta$ its vertex $L(\Delta)$ at which the function $f|_\Delta$
takes its supremum. Let $\nabla(\Delta)$ be the $N$-dimensional imaginary wedge
$\nabla_J$ in $\C^N\sm L$ with the corner at $L(\Delta)$, see the last
paragraph of \S \ref{salvet}. It is not difficult to show that the number of
bounded connected components $\Delta$ is equal to $\dim H_*(\C^N \sm L,
\Theta)$ (a proof of a more general statement see after Theorem \ref{reson2}
below). Statement C of Theorem \ref{reson1} is a corollary of the following
proposition.

\begin{proposition}[see \cite{VGZ}]
\label{basis} If all monodromy coefficients $\tau_j,$ $j=1, \dots, m,$ are
different from $1$, then

a) the group $H_N(\C^N \sm L, \Theta)$ is freely generated by the classes of
``double loops'' $\kappa(\Delta)$ corresponding to all bounded components
$\Delta$;

b) the group $H_N^{lf}(\C^N \sm L, \Theta)$ is freely generated by the classes
of imaginary wedges $\nabla(\Delta)$ corresponding to these components.
\end{proposition}

We shall prove simultaneously assertion a) for our local system $\Theta$ and
assertion b) for the dual system $\Theta^*$. Namely, let us lift the wedges
$\nabla(\Delta)$ to the total space $M^*$ of the latter system. Then the
intersection indices $\langle \kappa(\Delta_1), \nabla(\Delta_2)\rangle$ are
well-defined for all pairs of components $\Delta_1, \Delta_2$. It is easy to
calculate that the matrix consisting of all such intersection indices is
triangular with nonzero numbers on the diagonal; in particular it is
non-degenerate. \quad $\square$
\medskip

Many assertions of Theorem \ref{reson1} and Proposition \ref{basis} can be
extended to more general situations. However, it this case the set of resonant
values $T=(\tau_1, \dots, \tau_m)$ becomes more complicated.

Namely, let us consider the variety $\hat L \subset \CP^N$ consisting of the
support $L \subset \C^N$ of our arrangement and the improper plane. If this
variety is not a divisor with normal crossings then we take some its resolution
$(W,\tilde L) \to (\CP^N,\hat L)$ in the sense of \cite{hiro}. Spaces $W\sm
\tilde L$, $\CP^N \sm \hat L$ and $\C^N \sm L$ are diffeomorphic, hence we can
lift the local system $\Theta$ to $W \sm \tilde L$. For any irreducible
component $\tilde L_j$ of $\tilde L$, the monodromy coefficient corresponding
to a small circle around it is a monomial of the initial coefficients $\tau_1,
\dots, \tau_m$. Our collection $T=(\tau_1, \dots, \tau_m)$ is called {\em
resonant} if the value of some of these monomials is equal to 1 for {\em any}
resolution of $\hat L$. The trivial monomials $\tau_1, \dots, \tau_m$ and
$\tau_0^{-1} \equiv (\tau_1 \cdot \ldots \cdot \tau_m)^{-1}$ always appear
among these monomials: they correspond to proper images of initial planes $L_j$
and the improper plane. Proposition \ref{gen} implies immediately that if the
set $T$ is non-resonant then the map (\ref{canon}) is an isomorphism.

\begin{theorem}[see \cite{VGZ}]
\label{reson2} Suppose that our hyperplane arrangement $L \subset \C^N$ is the
complexification of a real one and has only normal crossings in $\C^N$, and all
numbers $\tau_1, \dots, \tau_m$ are not equal to 1. Then

A. Both groups $H_i(\C^N \sm L,\Theta)$, $H_i^{lf}(\C^N \sm L,\Theta)$ are
nontrivial only for $i=N$ and their dimensions are equal to the number of
bounded components of $\R^N \sm L$.

B. The group $H_i(\C^N \sm L,\Theta)$ is freely generated by the classes of
double loops corresponding to all these components, and the group
$H_i^{lf}(\C^N \sm L,\Theta)$ is freely generated by the classes of
corresponding complex wedges $\nabla(\Delta)$.

C. If moreover the collection of numbers $T=(\tau_1, \dots, \tau_m)$ is
non-resonant, then the group $H_i^{lf}(\C^N \sm L,\Theta)$ is freely generated
also by the classes of these bounded components.
\end{theorem}

A proof of the statement A. follows from the cell decomposition considered in
the last part of \S \ref{salvet}. Indeed, the corresponding cellular complex
calculating the homology group $H_i^{lf}(\C^N \sm L,\Theta)$ becomes, after
some norming of its canonical generators, isomorphic to the order complex of
our arrangement $\LL$ or, equivalently, of its real part $\LL \cap \R^N \equiv
\{L_j \cap \R^N\}$. By formula (\ref{abs}), the latter complex is homotopy
equivalent to the support $L\cap \R^N$ itself, which is obviously homotopy
equivalent to the bouquet of spheres corresponding to all these bounded
components. The proof of statements B. and C. is the same as for the case of
generic arrangements.
\medskip

\noindent {\bf Remark.} An analog of the deRham theory calculating the homology
of complex manifolds (e.g. of spaces $\C^N \sm L$) with coefficients in local
systems was developed in \cite{Deligne 70}, see also \cite{EShV}, \cite{Aomoto
77} and \S 8 in \cite{Yuz3}.

\section{Matroids and configuration spaces}
\label{matrs}

The simplest example of a configuration space is the space $M(N,2)$ of
collections of $N$ different points in $\C^1,$ see \S \ref{examples}. In more
complicated examples, we can consider collections of subvarieties in some
manifolds; a configuration space is the family of all collections such that the
pair (the manifold, the union of varieties) has a fixed topological type. Such
configuration spaces appear often in integral geometry and theory of special
functions, see \cite{Pham 67}, \cite{Gelfand 86}, \cite{VGZ}, \cite{ramif}. An
important class of such functions, called {\em general hypergeometric
functions}, appears if all our varieties are hyperplanes in $\C^N$ or $\CP^N$,
see \cite{Gelfand 86}, \cite{varchenko}. In this case the topology of
configuration spaces has especially deep relations with the algebraic geometry.
It is convenient to formulate these relations in terms of the theory of {\em
matroids}.

The notion of a matroid is a formalization of dimensional properties of a
central hyperplane arrangement, see \cite{maclane} and \cite{GS}.
\medskip

\noindent {\bf Definition.} A {\em matroid} is a finite set $U$ supplied with a
natural-valued function $r$ on the set $2^U \sm \{\emptyset\}$ of all non-empty
subsets of $U$ such that

1) for any such subset $I \in U$ we have $1 \le r(U) \le $ (the cardinality of
$I$);

2) if $I \subset J$ then $r(I) \le r(J)$;

3) for any $I, J$ we have $r(I \cap J) + r(I \cup J) \le r(I)+ r(J).$
\medskip

Any central hyperplane arrangement $\LL = \{L_1, \dots, L_m\}$ defines a
matroid, whose elements correspond to the planes $L_i$ and for any collection
$I$ of these elements $r(I)$ is the codimension of the intersection $L_I$ of
all planes from this collection.

In this case the arrangement $\LL$ is called {\em a realization} of the
corresponding matroid.

There exist matroids having {\em complex} realizations (i.e. realizations by
collections of complex hyperplanes in $\C^N$) but no real realizations, and
matroids having realizations over finite fields but no complex realizations,
etc. (On the other hand, any realization of a matroid over some field defines
also its realization over any extension of this field: in particular the
complexification of a real realization provides a complex realization).

Specifical properties of real affine hyperplane arrangements (roughly speaking,
the fact that two 0-dimensional planes $L_I$, $L_{I'}$ can lie either to one
side of any hyperplane not connecting them or to different sides) are
formalized in the notion of an {\em oriented matroid}, see \cite{ormat}.

The next important question is as follows: given a matroid, what is the set of
all its realizations over a given field $F$?

The study of realizations by central hyperplane arrangements in $\C^3$ (or,
equivalently, of arbitrary line arrangements in $\CP^2$) is deeply connected
with the theory of integer algebraic varieties (i.e. complex varieties defined
by equations with integer coefficients), see \cite{mnev}, \cite{mnev2}.

Accordingly to \cite{mnev}, \cite{mnev2}, for any integer algebraic subvariety
in some $\C^n$ there exists a matroid such that the space of its realizations
by plane arrangements in $\C^3$ is homotopy equivalent to our subvariety.

This relation allows one to construct spaces of realizations having very
delicious properties: these properties reflect the similar properties of the
corresponding algebraic varieties.

For instance, the equation $x^2 = -1$ is related with a matroid having complex
realizations but no real realizations. The equation $x^2 + y^2 = 0$ is related
with such a matroid that the real dimension of the space of its real
realizations is less than the complex dimension of the space of its complex
realizations; in particular the latter space is singular, see \cite{VS}.

The first of these examples (corresponding to the equation $x^2=-1$) is
constructed as follows.

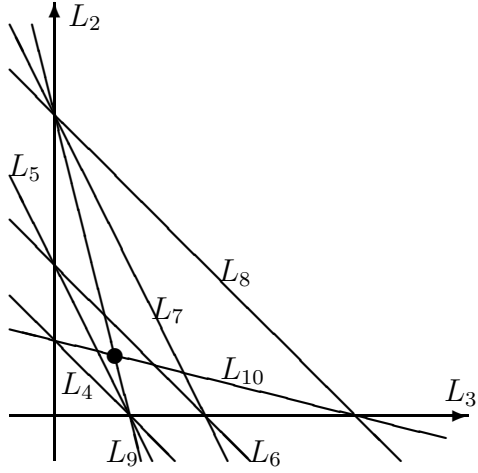
\begin{figure}
\begin{center}
\unitlength 1.00mm
\begin{picture}(61.00,61.00)
\put(0.00,6.00){\vector(1,0){61.00}} \put(6.00,0.00){\vector(0,1){61.00}}
\put(0.00,22.00){\line(1,-1){22.00}} \put(19.00,0.00){\line(-1,2){19.00}}
\put(0.00,32.00){\line(1,-1){32.00}} \put(29.00,0.00){\line(-1,2){29.00}}
\put(0.00,52.00){\line(1,-1){52.00}} \put(58.00,3.00){\line(-4,1){58.00}}
\put(3.00,58.00){\line(1,-4){14.50}}
\put(10.00,59.00){\makebox(0,0)[cc]{$L_2$}}
\put(60.00,9.00){\makebox(0,0)[cc]{$L_3$}}
\put(9.00,10.00){\makebox(0,0)[cc]{$L_4$}}
\put(2.00,39.00){\makebox(0,0)[cc]{$L_5$}}
\put(34.00,1.00){\makebox(0,0)[cc]{$L_6$}}
\put(21.00,20.00){\makebox(0,0)[cc]{$L_7$}}
\put(30.00,25.00){\makebox(0,0)[cc]{$L_8$}}
\put(15.00,1.00){\makebox(0,0)[cc]{$L_9$}}
\put(31.00,12.00){\makebox(0,0)[cc]{$L_{10}$}}
\put(14.00,14.00){\circle*{2.00}}
\end{picture}
\end{center}
\caption{Imaginary line arrangement} \label{mne}
\end{figure}

The first four elements $L_1, \dots, L_4$ are {\em in general position} (i.e.,
$r(\{I\})$ is equal to 1 for any 1-element set $I \subset \{1,2,3,4\}$, to 2
for any 2-element set and to 3 for 3- and 4- element sets). In the terms of
possible realizations by lines in $\CP^2$ or $\RP^2$ this means that these
lines meet in general position in the standard sense: none three of them meet
at one point. Given a realization of our matroid including these four elements,
they fix a coordinate system in $\PP^2$: we can take $L_1$ as the line ``at
infinity'', $L_2$ as the line $\{x=0\}$, $L_3$ as the line $\{y=0\}$, and
choose the scaling of coordinates in such a way that the points $L_2 \cap L_4$
and $L_3 \cap L_4$ have coordinates $(0,1)$ and $(1,0)$ respectively.

Further we add the element $L_5$ with unique non-generic condition $r(L_3, L_4,
L_5)=2$, i.e. the corresponding three lines of any realization in $\PP^2$
should intersect at one point. The intersection of this line $L_5$ with $L_2$
will be some point with coordinates $(0,\alpha)$. The next line $L_6$ should
pass through this intersection point and be parallel to $L_4$ (in the terms of
the matroid these conditions are expressed as $r(L_2, L_5, L_6)=2$ and $r(L_1,
L_4, L_6)=2$ respectively). The point $L_6 \cap L_3$ will have coordinates
$(\alpha,0)$. The next line $L_7$ passes through this point and is parallel to
$L_5$, i.e. we have conditions $r(L_3,L_6,L_7)=2$ and $r(L_1,L_5, L_7)=2$. Its
intersection point with the line $L_2$ has coordinates $(0,\alpha^2)$. The next
line $L_8$ passes through the last point and is parallel to $L_6$ and $L_4$,
i.e. we have $r(L_2,L_7,L_8)=2$ and $r(L_1,L_4,L_6,L_8)=2$. Its intersection
point with the line $L_3$ has coordinates $(\alpha^2,0)$. The next line $L_9$
should pass through the points $(1,0)$ and $(0,\alpha^2)$, i.e.
$r(L_3,L_4,L_5,L_9)=2$ and $r(L_2,L_7,L_8,L_9)=2.$ The line $L_{10}$ should
pass through the points $(0,1)$ and $(\alpha^2,0)$, i.e. $r(L_2,L_4,L_{10})=2$
and $r(L_3,L_8,L_{10})=2$.

Finally, we claim that the lines $L_9$ and $L_{10}$ are parallel (i.e.
$r(L_1,L_9,L_{10})=2$): the intersection point covered by the ``black hole'' in
Fig.~\ref{mne} should lie at the infinity. This is possible only if
$\alpha^4=1$. But $\alpha \ne 1$ (because $L_5 \ne L_4$) and $\alpha^2 \ne 1$
(because $L_8 \ne L_4$). Therefore any complex realization of this matroid
corresponds to a number $\alpha$ with $\alpha^2=-1.$

\section{Applications in integral geometry:
general hypergeometric functions} \label{hgf}

The theory of {\em general} (or {\em multidimensional}) {\em hypergeometric
functions} was initiated by K.~Aomoto \cite{Aomoto 75}, \cite{Aomoto 77} and
I.M.~Gelfand \cite{Gelfand 86}. These functions form an important class of
functions given by integral transformations: they represent all these functions
to the same extent as the plane configurations represent all algebraic
varieties.

The starting point of this theory is the Gauss' hypergeometric integral
\begin{equation}
\label{0.gauss} \Gamma(a;\alpha_1,\alpha_2, \alpha_3) \equiv \int_0^1
z^{\alpha_1}(z-1)^{\alpha_2}(z-a)^{\alpha_3} dz,
\end{equation}
where $\alpha_j$ and $a$ are complex numbers, $a \not \in [0,1]$\footnote{Gauss
himself wrote it in a different but equivalent form}.

This integral converges absolutely if ${\rm Re}\, \alpha_1>-1,$ ${\rm Re}\,
\alpha_2>-1.$ For any fixed $a$ it is a holomorphic function in the domain of
the space $\C^3$ of exponents $\alpha= (\alpha_1, \alpha_2, \alpha_3)$
distinguished by these inequalities. The analytical continuation of this
function to the space of all $\alpha \in \C^3$ is a meromorphic single-valued
function whose poles are the hyperplanes on which either $\alpha_1$ or
$\alpha_2$ is a negative integer number.

On the other hand, fixing $\alpha$ and moving $a$ we obtain an analytical
function on $a$ with a ramification at points $0$ and $1$, satisfying the
famous hypergeometric equations, see \cite{BE}. If the exponents $\alpha$ are
generic, then the space of solutions of these equations at any point $a$ is
two-dimensional and is spanned by different leaves of the analytic continuation
of this function.

Much more generally, we can consider the integral
\begin{equation}
\label{0.gelf} \Gamma(\lambda;\alpha) = \int_\Delta f_1^{\alpha_1}\cdot \ldots
\cdot f_k^{\alpha_k} dz_1 \wedge \ldots \wedge dz_n,
\end{equation}
where $f_j$ are some polynomial functions $\C^N \to \C,$ depending analytically
on parameters $\lambda \in \C^m,$ $\alpha_j$ are complex exponents, and
integration is taken over some relatively closed (i.e. locally finite) but,
generally, not compact (i.e. not finite) $n$-dimensional cycles in the space of
non-zero values of the integration function (or, more precisely, in the space
of a covering over this space in which our function becomes single-valued).

An important class of such problems is as follows. We fix some matroid and
allow the functions $f_j$ be the linear functions $\C^N \to \C$ whose zero sets
form all possible complex realizations of this matroid.

The integral (\ref{0.gelf}) also defines a meromorphic function on $\alpha$
(for fixed $a$), cf. \cite{BG}, and a branching analytical function on
$\lambda$ for fixed $\alpha$; the set of its ramification is the set of such
values of the parameter $\lambda,$ that the corresponding functions $f_j$ have
the topologically ``non-generic'' sets of zeros (i.e. the rank function $r$ of
the plane arrangement ``jumps'').

Again, the main problems here are as follows.

1) to describe the polar set of these integrals, considered as meromorphic
functions on $\alpha$;

2) to describe the ramification of these integrals for fixed $\alpha$ and
moving parameter $\lambda$ of the set of functions $f_j$;

3) to find the number of linearly independent functions on the parameter spaces
$\C^m,$ given by such integrals.

The solution of these analytical problems is closely related with the following
topological ones.

A) the calculation of homology groups related with our collections of functions
$f_j$ and containing all possible integration chains $\Delta$;

B) the study of the maps similar to (\ref{canon}) for such groups;

C) the study of the topology (especially of the fundamental group) of the
configuration space of all realizations of our matroid;

D) the study of the {\em homology vector bundle} over this configuration space,
whose fiber over an arrangement is the corresponding homology group considered
in A); especially the study of the {\em monodromy representation} of the
fundamental group from C) in the fiber of this bundle.

Let us consider all these problems (and their applications) for the classical
integral (\ref{0.gauss}). In this case the solution of the problem 1) follows
from Fig.~\ref{doloop}a): if $\alpha_1$ and $\alpha_2$ are not integer, then
our improper integral along the interval $(0,1)$ can be replaced by a similar
integral along the ``double loop'' taken with the coefficient $1/(1-e^{2\pi i
\alpha_1})(1-e^{2\pi i \alpha_2})$, which is regular.

Similarly, in the general case of the integral (\ref{0.gelf}) if the collection
of exponents $\tau_j \equiv e^{2 \pi i \alpha_j}$ is not resonant in the sense
of \S \ref{twist} then any integration cycle can be {\em regularized} by a
``multidimensional double loop'', see \cite{ramif}.
\medskip

\noindent {\bf Remark.} The latter assertion is not a formal corollary of
Proposition \ref{gen}, i.e. of the invertibility of the map (\ref{canon}) for
non-resonant exponents. Indeed, in the study of integrals as functions of these
exponents, the integration cycles cannot be considered as elements of the group
\begin{equation}
H^{lf}_N(\C^N \sm L, \Theta), \label{125}
\end{equation}
because the definition of this group involves some factorization depending on
exponents. Still, the geometrical construction of ``double loops'' works even
in this case and allows us to regularize the integrals.
\medskip

On contrary, if we fix the exponents $\alpha_j$ then we can consider the
integration cycles as elements of the group (\ref{125}). The dimension of the
space of different integral functions (considered as functions on parameters
$\lambda$ of functions $f_j$ only) is then no greater than the dimension of
this group. This estimate can be not sharp in some exotic examples, however for
many important matroids it is sharp. For instance it is so for matroids
corresponding to generic arrangements and, moreover, to all arrangements with
normal crossings in $\C^N$ but generally not in $\CP^N$, see \cite{VGZ}.

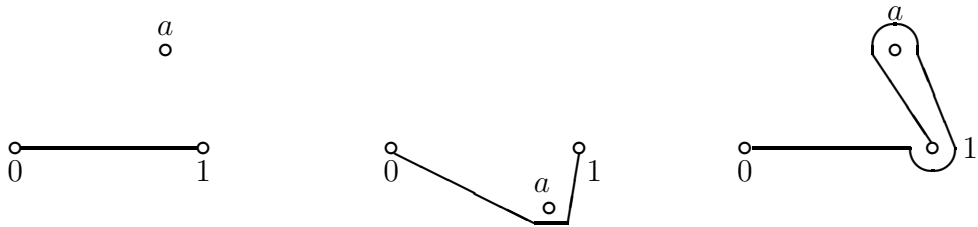
\begin{figure}
\begin{center}
\unitlength 1.00mm \special{em:linewidth 0.4pt}
\begin{picture}(127.00,31.00)
\put(0.75,13.00){\line(1,0){23.50}} \put(25.00,13.00){\circle{1.50}}
\put(0.00,13.00){\circle{1.50}} \put(75.00,13.00){\circle{1.50}}
\put(122.00,13.00){\circle{1.50}} \put(50.00,13.00){\circle{1.50}}
\put(97.00,13.00){\circle{1.50}} \put(0.00,10.00){\makebox(0,0)[cc]{$0$}}
\put(25.00,10.00){\makebox(0,0)[cc]{$1$}}
\put(50.00,10.00){\makebox(0,0)[cc]{$0$}}
\put(97.00,10.00){\makebox(0,0)[cc]{$0$}}
\put(77.00,10.00){\makebox(0,0)[cc]{$1$}}
\put(127.00,13.00){\makebox(0,0)[cc]{$1$}} \put(20.00,26.00){\circle{1.50}}
\put(117.00,26.00){\circle{1.50}} \put(71.00,5.00){\circle{1.50}}
\put(117.00,31.00){\makebox(0,0)[cc]{$a$}}
\put(70.00,8.00){\makebox(0,0)[cc]{$a$}}
\put(20.00,29.00){\makebox(0,0)[cc]{$a$}} \put(50.50,12.25){\line(2,-1){18.50}}
\put(69.00,3.00){\line(1,0){4.46}} \put(75.00,12.25){\line(-1,-6){1.54}}
\put(98.00,13.00){\line(1,0){21.00}} \put(122.00,13.00){\oval(6.00,6.00)[b]}
\put(117.00,25.50){\oval(6.00,8.00)[t]} \put(125.00,13.00){\line(-2,5){5.00}}
\put(114.00,25.50){\line(2,-3){7.80}}
\end{picture}
\end{center}
\caption{Monodromy of the hypergeometric integral} \label{mon}
\end{figure}

The proof (see \cite{ramif}) is based on the study of the monodromy action in
the homology bundle: starting from a single integration cycle and acting on it
by all elements of the monodromy group we can obtain a collection of cycles
generating the whole group (\ref{125}).
\medskip

\noindent {\bf Example.} In the case of the integral (\ref{0.gauss}) the
configuration space is the punctured plane $\C^1 \sm \{0,1\}$ of all admissible
values of the parameter $a$. Let us move this value along a closed loop in this
space going around some singular point, say the point $1$. Deforming
simultaneously the integration cycle $(0,1)$ in such a way that at no instant
it intersects singular points $0,1$ and the current value of $a$, we get the
cycle shown in Fig.~\ref{mon} right. The integral of type (\ref{0.gauss}) along
this cycle is equal to the analytical continuation of the initial integral
(\ref{0.gauss}) along our loop in the parameter space. The difference between
this cycle and the initial one is equal to the interval $(1,a)$ passed in two
opposite directions {\em on two different leafs} of the Riemann surface on
which our integration form is single-valued. If $\alpha_3$ is not integer, then
the integral along this cycle is not identically equal to zero, and we get the
second integral function independent on the first one. Since the dimension of
the group $H_1^{lf}(\C^1 \sm \{0,1,a\},\Theta)$ is equal to 2, the lower and
upper estimates on the number of linearly independent integral functions
coincide and problem 3) is solved in this particular case.

On algebraic properties of analytic functions of this kind see, in particular,
\cite{DM 86} and \cite{GKZ}.

Similar (although much more complicated) methods allow us to prove analogous
results in multidimensional situations, see \cite{VGZ}, \cite{ramif}. In
particular if the sets $\{x|f_j(x)=0\}$ form hyperplane arrangements with
normal crossings in $\C^N$ and the set of exponents $\tau_j \equiv e^{2 \pi i
\alpha_j}$ is non-resonant, then the homological estimate also is sharp, i.e.
we have the full number of linearly independent integral functions.

These and similar functions have wonderful applications in mathematical
physics, see e.g. \cite{Schechtman-Varchenko 91}, \cite{varchenko}.

\section{How if the collection of planes is infinite ?}

There are many important topological subspaces in $\R^N$ that can be
represented as continuous families of planes. A great source of such spaces is
the {\em discriminant theory}, see \cite{b94}, \cite{fasis}. The
above-described strategy of the topological investigation of unions of planes
(and complements of such unions), based on simplicial resolutions, works also
in these cases (after appropriate modification). In particular we need to
consider {\em continuous order complexes} and construct the {\em conical
resolutions} of such spaces. Here we describe a simple example: that of
determinant varieties.

Let $\KK$ be any of fields $\R, \C$ or $\HH.$ The {\em determinant} variety
$Det(\KK^n) \subset End(\KK^n) \sim \KK^{n^2}$ consists of all degenerate
operators $\KK^n \to \KK^n.$

We construct a resolution of this variety that provides a calculation of its
Borel--Moore homology groups and, by the Alexander duality, the {\em most
complicated calculation} of cohomology groups of their complementary spaces
$GL(\KK,n)$.

The {\em tautological resolution} of $Det(\KK^n)$ is defined by elimination of
quantifiers (which is an analog of ``taking the sets $S_j$ separately'' in the
justification of the inclusion--exclusion formula). Namely, an operator $A$
belongs to $Det(\KK^n)$ if $\exists$ a point $x \in \KP^{n-1}$ such that $\{x\}
\subset \ker A$. Define the resolution space $det_1(\KK^n)$ as the space of all
pairs $(x, A) \in \KP^{n-1} \times End(\KK^n)$ such that $\{x\} \in \ker A$.
This space admits the (tautological) structure of a $(n^2-n)$-dimensional
$\KK$-vector bundle over $\KP^{n-1}$, whose fiber $L(x)$ consists of all $A$
such that $\{x\} \in \ker A$. The obvious projection $\pi: det_1(\KK^n) \to
Det(\KK^n)$ is regular over operators with one-dimensional kernels, but the
pre-image of an operator with $\dim \ker=l$ is isomorphic to $\KP^{l-1}$.

The situation is very similar to the one considered in the arrangement theory:
the variety $Det(\KK,n)$ is the union of spaces $L(x)$ in the same way as the
space $L$ was the union of planes $L_i$. Keeping the analogy, we construct the
order complex of all intersections of these spaces $L(x)$. It is not
straightforward because the family of planes $L(x)$ is not discrete, and
moreover the set of such planes passing through one and the same point of
$Det(\KK^n)$ can be continuous. Indeed, the possible intersections of several
planes $L(x_j) \subset End(\KK^n)$ are just the planes of the form $L(X)$ where
$X$ is a subspace of a certain dimension in $\KK^n$ (i.e. a point of a certain
Grassmannian manifold $G_i(\KK^n),$ $i\in [1,n]),$ and $L(X)$ consists of all
operators whose kernels contain $X$. \medskip

Thus our poset of all planes and their intersections is the {\em disjoint union
of all Grassmann manifolds} $G_1(\KK^n), \ldots, G_{n-1}(\KK^n), G_n(\KK^n).$
The {\em continuous order complex} of all these Grassmannians is defined as
follows. Consider the join $G_1(\KK^n)* \ldots * G_{n}(\KK^n),$ i.e., roughly
speaking, the naturally topologized union of all simplices whose vertices
correspond to points of different Grassmannians. Such a simplex is {\em
coherent} if the planes corresponding to its vertices form a flag. The desired
order complex $\Xi(\KK^n)$ is the union of all coherent simplices, with
topology induced from that of the join. This is a cone with vertex $\{\KK^n\}
\in G_n(\KK^n).$ Its {\em link} $\partial \Xi(\KK^n)$ is the union of coherent
simplices not containing the vertex $\{\KK^n\}$.

This link $\partial \Xi (\KK^n)$ is homeomorphic to the sphere $S^M,$
$M=\frac{1}{2}n(n-1)(\dim_{\R}\KK) + n-2$. (Probably this fact is assumed in
Remark 1.4 of \cite{bs}, see also \cite{V-2}, \cite{fasis}.) Hence $\Xi(\KK^n)$
is homeomorphic to a ball.

The {\em conical resolution} of $Det(\KK^n)$ is constructed as a subset of the
direct product $\Xi(\KK^n) \times Det(\KK^n).$ To any plane $X \subset \KK^n$
there corresponds a subspace $\Xi(X) \subset \Xi(\KK^n),$ namely, the union of
all coherent simplices all whose vertices correspond to planes lying in $X$.
This is a cone with vertex $\{X\}$, and is homeomorphic to a closed ball.
Define the conical resolution $\delta(\KK^n) \subset \Xi(\KK^n) \times
Det(\KK^n)$ as the union of the products $\Xi(X) \times L(X)$ over all planes
$X$ of dimensions 1, \ldots, $n$. It is easy to see that the obvious projection
$\delta(\KK^n) \to Det(\KK^n)$ induces a homotopy equivalence of one-point
compactifications of these spaces (indeed, this projection is proper and
semialgebraic, and all its fibers are contractible cones of the form $\Xi(X)$).
On the other hand, the space $\delta(\KK^n)$ has a nice filtration: its term
$F_i$ is the union of products $\Xi(X) \times L(X)$ over planes $X$ of
dimensions $\le i$. The term $F_i \setminus F_{i-1}$ of this filtration is the
total space of a fibre bundle over $G_i(\KK^n).$ Its fiber over a point $\{X\}$
is the space $(\Xi(X) \setminus \partial \Xi(X)) \times L(X)$, and is
homeomorphic to an Euclidean space. Thus the Borel--Moore homology group of
this term can be reduced to that of the base. The spectral sequence, generated
by this filtration and converging to the Borel--Moore homology group of
$Det(\KK^n)$ (or, equivalently, to the cohomology group of the complementary
space $GL(\KK^n)$), degenerates at the first term (i.e. $E_1^{p,q} \equiv
E_\infty^{p,q}$) and gives, in particular, the homological {\em Miller
splitting}
\begin{equation}
\label{miller} H_m(GL(\C^n)) = \bigoplus_{k=0}^n H_{m-k^2}(G_k(\C^n))
\end{equation}
and similar splittings for $\KK= \R$ and $\HH$.

There is a plenty of other problems in which the technology of conical
resolutions works, see \cite{kotor}, \cite{siersma}. Among them are the theory
of knots and generic plane curves (see the next section), topological study of
spaces of continuous maps, of smooth functions without complicated
singularities, of operators with simple spectra, of nonsingular projective
hypersurfaces...

\section{Applications and analogies in differential topology}

The space $M(N,2)$ (see \S \ref{examples}) can be considered as the space of
all embeddings to $\C^1$ of a finite set of cardinality $N$. In a similar way
we can consider the space of all smooth embeddings $S^1 \hookrightarrow \R^n$,
i.e. regular knots in $\R^n$; if $n=3$ then the 0-dimensional cohomology
classes of this space are the knot invariants. We can study these and other
cohomology classes (in the case of any $n\ge 3$) by essentially the same
methods as in \S \ref{sires}, \ref{sthoty} (curiously, it was done earlier, see
\cite{knsp}). Consider the space $\K$ of all smooth maps $S^1 \to \R^n,$ define
the discriminant $\Sigma \subset \K$ as the space of all maps that are not
smooth embeddings, and study the group $H^*(\K \sm \Sigma)$. To do it, we take
a conical resolution of the discriminant set $\Sigma$. It is possible because
this set is swept out by a reasonable family of subspaces in $\K$. These
subspaces are parameterized by all unordered pairs of points $(x,y) \subset
S^1.$ By obvious reasons these pairs are called {\em chords}, they run over the
2-dimensional {\em chord space} $\overline{B(S^1,2)}.$

For any such pair the corresponding subspace $L(x,y)$ consists of all maps
$f:S^1 \to \R^n$ such that $f(x)=f(y)$ if $x \ne y$ or $f'(x)=0$ if $x=y$. Such
subspaces form the tautological resolution of $\Sigma$. Then we take the order
complex of all possible intersections
\begin{equation}
\label{peres} L(x_1, y_1) \cap L(x_2,y_2) \cap ...
\end{equation}
and limit positions of such intersections (all of them are subspaces in $\K$
whose codimensions are multiples of $n$), {\em supply it with a natural
topology}, and define the conical resolution in exactly the same way as
previously, i.e. as a subset of the direct product of this order complex and
the space $\K$. Then we define the filtration on this resolution by the
codimensions (divided by $n$) of these planes and consider the arising spectral
sequence.

The homological study of the arising resolution space is known as the theory of
{\em finite-type knot invariants}, see \cite{bnbib}, and different its
generalizations, including the (equally interesting) calculation of higher
dimensional cohomology classes of spaces of knots, see \cite{knsp},
\cite{royal}, \cite{turdis}. In Table 1 we give a short list of parallel
notions and objects in both theories. Of course, a large part of its right half
cannot be explained here; this table is rather a kind of Rosetta stone for
those combinatorialists who will study the theory of knot invariants, see
\cite{bnbib}.

\begin{table}
\begin{tabular}{|l|l|}
\hline
Theory of arrangements & Knot theory \\
\hline
Space $\R^N$ & Space ${\mathcal K} = C^\infty(S^1,\R^n)$ \\
\hline Union of planes $L=\cup L_i \subset \R^N$ & Discriminant subset
$\Sigma \subset {\mathcal K}$ \\
\hline
Set of indices $\{1, \ldots, m\}$ & Chord space $\overline{B(S^1,2)}$ \\
\hline
A plane $L_i$ & A subspace $L(x,y)$, $x,y \in S^1$ \\
\hline Disjoint union of planes $L_i$ & Tautological resolution $F_1 \sigma$
of $\Sigma$ \\
\hline Simplicial resolution $L'$ of $L$ & Conical resolution
$\sigma$ of $\Sigma$ \\
\hline Subsets $I \subset \{1, \ldots, m\}$ & Combinatorial types of chord
\\
with codim$L_I =p$ & configurations $J$ with \\ & codim$L(J) = pn$ \\
\hline
A prism $L'_I$ & A $J$-block in $\sigma$ \\
\hline K\"unneth isomorphism for &
Thom isomorphism for the fibration \\
homology of $\check L'_I = \check \Delta(I) \times L_I$ &
of pure $J$-blocks by spaces $L(J')$ \\
\hline Shuffle product formulas of \S \ref{mumu} & Multiplication formulas
by   \\
& Kontsevich (for invariants)  \\ & and Turchin (for higher \\ & cohomology classes) \\
\hline Integral representations of
\S\S \ref{basbas}, \ref{osos} & Kontsevich integral \\
\hline
Homotopy splitting (\ref{rel}) & Kontsevich stabilization theorem \\
\hline Combinatorial realization of \S \ref{cococo} & Combinatorial formulas
for  \\
& invariants \cite{PV}, \cite{GPV} and other \\ & cohomology classes \cite{vtt} \\
\hline
\end{tabular}
\caption{The analogy between the arrangement theory and the knot theory}
\label{glos}
\end{table}

It is necessary also to mention the exceptional value of the Arnold's identity
(\ref{relation}) for the construction of the Kontsevich's integral
\cite{konec}, \cite{Konn}.

Of course, the space ${\mathcal K}$ is infinite-dimensional, and formally we
cannot use the Alexander duality in it: the usual (i.e. finite-dimensional)
cohomology classes of the space of knots ${\mathcal K} \setminus \Sigma$ should
correspond to ``infinite-dimensional cycles'' in $\Sigma$, whose definition
requires some effort. The strict construction of such cycles corresponding to
finite-type cohomology classes uses the techniques of finite-dimensional
approximations, see \cite{knsp}.

Similarly, we can consider the space of smooth embeddings of finitely many
circles into $\R^n$, it gives us the theory of finite-type cohomology classes
of spaces of links.

One more space of this type is that of all plane curves without triple points,
see \cite{Arnold 94}, \cite{Arnold 95}, \cite{congr}, \cite{ornam},
\cite{merx}--\cite{M}, \cite{shum}, \cite{hov}, \cite{modov}. It is related
very much with the arrangements $A(N,3)$, see \cite{BW} and \S \ref{examples}.
\bigskip

I thank all people to whom I owe my knowledge of the described area, including
V.I. Arnold, A. Bj\"orner, V. Welker, A.M. Vershik, I.M. Gelfand, R.T.
Zhivaljevich, A.V. Zelevinskii, M.M. Kapranov, A.B. Merkov, H.E. Mnev, S.A.
Yuzvinsky, G.L. Rybnikov, V.V. Serganova, V.E. Turchin, G.M. Ziegler, B.Z.
Shapiro, and K. Schultz.


\begin{thebibliography}{99}

\bibitem{Aomoto 75} Aomoto K. Les \'equations aux diff\`erences
lin\`eaires des fonctions multiformes // J.~Fac.~Sci.~Univ. of Tokyo 22 (1975),
271--297.

\bibitem{Aomoto 77} Aomoto K. On the structure of integrals of
power products of linear functions // Sci. papers, Coll. Gen. Ed., Univ. Tokyo.
--- 1977. V.~27, p.~49--61.

\bibitem{ararr} Arnold V.I. The cohomology ring of the group
of colored braids.
Math. Notes {\bf 5} 138--140.

\bibitem{Arnold-4} Arnold V.I. On some topological invariants of the
algebraic functions. Trans. Moscow Math. Soc. {\bf 21} (1970), 30--52.

\bibitem{Arnold 94} Arnold V.I. Plane curves, their invariants,
perestroikas and classifications, in V.~I.~Arnold (ed.), {\em Singularities and
Bifurcations,} Adv. in Soviet Math. {\bf 21}, AMS, Providence RI, 1994, 33--91.

\bibitem{Arnold 95} Arnold V.I.
Invariants and perestroikas of fronts in the plane, Proc. Steklov Math. Inst.
{\bf 209}, 14--64.

\bibitem{AVGL} Arnold V.I., Vassiliev V.A., Goryunov V.V.,
and Lyashko O.V. Singularities II. Applications. Itogi Nauki i Tekhn. VINITI.
Fundam. napravl. Moscow, VINITI, 1989. English transl.: Encycl. Math. Sci.,
vol. 6, Springer-Verlag, Berlin and New York, 1993.

\bibitem{Artin} Artin E. Theorie der Z\"opfe.// Abh. Math. Semin. Univ.
Hamburg. 1925, 4, S.~47--72. (English Transl.: Theory of braids. Ann. Math.,
1947, 48:1, p.~101--126).

\bibitem{BBLSW} Babson E., Bj\"orner A., Linusson S., Shareshian J.,
and Welker V. The complexes of not $i$-connected graphs// Topology, {\bf 38}:2,
1999, 271--299.

\bibitem{bnbib} Bar-Natan, D. (1994--) {\em Bibliography
of Vassiliev Invariants}. Web publication {\tt
http://www.ma.huji.ac.il/\~{}drorbn/VasBib/VasBib.html}

\bibitem{BE} Higher Transcendental
Functions. (H.~Bateman project) Vol.~1. McGraw-Hill, New York-Toronto-London,
1953.

\bibitem{BG} Bernstein I.N. and Gelfand S.I., Meromorphy of the function
$P^\lambda$// Funct. Anal. and its Appl., 3, 68--69 (1969).

\bibitem{Bourbaki} Bourbaki N. Groupes et alg\`ebres de Lie,
Chapitres 4, 5, 6, Hermann, Paris, 1968.

\bibitem{ormat} Bj\"orner A., Las Vergnas M., Sturmfels B., White N.,
Ziegler G.M. Oriented matroids. Cambridge Univ. Press, London and New York,
1992.

\bibitem{bjorner} Bj\"orner A. Subspace arrangements//
In: Proc. 1-st European congr. of Mathematicians (Paris, 1992). Birkh\"auser,
1993. 321--370.

\bibitem{bjornhb} Bj\"orner A. Topological methods//
In: Handbook of Combinatorics. R.~Graham, M.~Gr\"otschel and L.~Lovasz, eds.
Amsterdam: North--Holland, 1995; p. 1819--1872.

\bibitem{Bjorner-2} Bj\"orner A. Nonpure shellability, $f$-vectors,
subspace arrangements and complexity// In: Formal power series and algebraic
combinatorics (New Brunswick, N.~J., 1994). Providence, RI, 1996. 25--53.
(DIMACS Ser. Discrete Math. Theoret. Comput. Sci., {\bf 24}.)

\bibitem{BLY} Bj\"orner A., Lov\'asz L. and Yao A. Linear
decision trees: volume estimates and topological bounds// In: Proc. 24-th Ann.
ACM Symp. on Theory of Computing, New York: ACM Press, 1992. 170--177.

\bibitem{BW} Bj\"orner A., Welker V., The
homology of ``$k$-equal'' manifolds and related partition lattices// Adv.
Math., 110:2 (1995), 277--313.

\bibitem{BZ} Bj\"orner A., Ziegler G.M., Combinatorial Stratification
of complex arrangements. Journal of the AMS. 1992. V.~5. N~1. P.~105--149.

\bibitem{bs} Borel A., Serre J.-P.
Cohomologie d'immeubles et de groupes $S$-arithmetiques,
// Topology {\bf 15} (1976), 211--231.

\bibitem{Brieskorn} Brieskorn E. Sur les groupes de tresses
(d'apr\`es V.~I.~Arnol'd)// In: Seminaire Bourbaki, 1971/72, No.~401. Berlin:
Springer, 1973. 21--44. (Lect. Notes Math, 317).

\bibitem{cart} Cartier P. Les arrangements d'hyperplans: un chapitre
de g\'eometrie combinatoire// in: Seminaire Bourbaki 1980/81", Lecture Notes in
Math. {\bf 901}, Springer Verlag, Berlin-Heidelberg-New York, 1981, 1--22.

\bibitem{CP} De Concini C., Procesi, C. Wonderful models
of subspace arrangements// Selecta Math., New Ser., {\bf 1}:3 (1995), 459--494.

\bibitem{DGM} Deligne P., Goresky M., MacPherson R.
L'algebre de cohomologie du compl\'ement, dans un espace affine, d'une famille
finie de sous-espaces affines// Michigan J. Math. {\bf 48} (2000)

\bibitem{Deligne 70} Deligne P., Equations diff\'erentielles
\`a points singuliers r\'eguliers// Lect. Notes Math. 163, Springer, Berlin,
1970.

\bibitem{DM 86} Deligne P., Mostow G., Monodromy of
hypergeometric functions and nonlattice integral monodromy// Publ. Math. IHES
63 (1986), 90 p.

\bibitem{Dimm} Deligne P., Les immeubles des groupes de tresses
g\'en\'eralis\'es// Invent. Math., 17 (1972), 273--302.

\bibitem{dLS} de Longueville M., Schultz, C. The cohomology
rings of complements of subspace arrangements, Math. Annalen (to appear)//
Preprint version available at: {\tt
http://www.math.fu-berlin.de/$\sim$cschultz/arrangements.html}

\bibitem{EShV} Esnault H., Shechtman V., Viehweg E.,
Cohomology of local systems on the complement of hyperplanes// Invent. Math.,
{\bf 109:}3 (1992), 557--561 and {\bf 112:2} (1993), 447.

\bibitem{FaN} Fadell E., Neuwirth L., Configuration spaces//
Math. Scand., 1962, 10:1, p.~111--118.

\bibitem{Folkman} Folkman J. The homology groups of a
lattice// J. Math. Mech. 1966 {\bf 15}, 631--636.

\bibitem{Fuchs} Fuchs D.B., Cohomology of the braid group mod 2 //
Funct. anal. and its Appl. 4:2 (1970), 143--152.

\bibitem{Gelfand 86} Gelfand I.M.,
A general theory of hypergeometric functions, Sov. Math. Doklady, 33, 1986,
573--577.

\bibitem{Gelfand-Gelfand 86} Gelfand I.M., Gelfand S.I.,
Generalized hypergeometric equations // Sov. Math. Doklady, 33, 1986.

\bibitem{GZ 86} Gelfand I.M., Zelevinskii, A.V., Algebraic and
combinatorial aspects of the general theory of hypergeometric functions
//
Funct. Anal. and its Appl., 20:3, 1986, 183--197.

\bibitem{GR} Gelfand I.M., Rybnikov G.L.,
Algebraic and topological invariants of oriented matroids
// Soviet Math. Doklady {\bf 40} (1990), 148--152.

\bibitem{GKZ} Gelfand I.M., Kapranov M.M., Zelevinskii A.V.,
Discriminants, Resultants, and Multidimensional Determinants. Birkhauser,
Boston MA, 1994, 523 p.

\bibitem{GS} Gelfand I.M., Serganova V.V., Combinatorial
geometries and torus strata on homogeneous compact manifolds
// Russian Math. Surveys {\bf 42} (1987), 133--168.

\bibitem{Gor} Goryunov V.V., Cohomology of braid groups of series
$C$ and $D$// Trans. Moscow Math. Soc., 2, 233--241 (1982).

\bibitem{GM} Goresky M., MacPherson R.,
Stratified Morse Theory, Springer, 1988, Berlin a.o.

\bibitem{GPV} Goussarov M., Polyak M., Viro O. Finite type
invariants of classical and virtual knots// Topology {\bf 39}:5, (2000)
1045--1068.

\bibitem{Hattori} Hattori A., Topology of ${\bf C}^n$ minus
a finite number of affine hyperplanes in general position// J.~Fac.~Sci.~Univ.
of Tokyo 22 (1975), 205--219.

\bibitem{hiro} Hironaka H.
Resolution of singularities of an algebraic variety over a field of
characteristic zero// Ann Math. 1964 {\bf 79}:1, 109--203 and {\bf 79}:2,
205--326.

\bibitem{JOS} Jewell K., Orlik P., Shapiro B.Z.,
On the complement of affine subspace arrangements// Topology and its
Applications {\bf 56} (1994), 215--233.

\bibitem{hov} Khovanov M., Doodle groups//
Trans. AMS {\bf 349} (1996), 2297--2315.

\bibitem{hov2} Khovanov M., Real $K(\pi,1)$ arrangements
from finite root systems// Mathematical Research Letters {\bf 3}, (1996),
261--274.

\bibitem{konec} Kontsevich M. Vassiliev's knot invariants//
In: Adv. in Sov. Math., {\bf 16:2} (1993), AMS, Providence RI, 137--150.

\bibitem{Konn} Kontsevich M., Formal (non-)commutative symplectic
geometry// In: L.~Corvin, I.~Gel'fand, J.~Lepovsky (eds.), The I.~M.~Gel'fand's
mathematical seminars 1990--1992, Birkh\"auser, Basel, 1993, 173--187.

\bibitem{maclane} Mac Lane S., Some interpretations of abstract
linear dependence in terms of projective geometry// Amer. J. Math., {\bf 58}
(1936), 236--240.

\bibitem{Lin} Lin V.Ya., Artin braids and related groups and spaces
// In: Itogi nauki i tekhn. Ser. Algebra, Geometry, Topology.
V. 17, Moscow.: VINITI, 1979, 159--227. English transl.: J. Soviet Math., 18:5
(1982).

\bibitem{merx} Merkov A.B. Finite order invariants of
ornaments// J. of Math. Sciences, {\bf 90}:4 (1998), 2215--2273.

\bibitem{merxx} Merkov A.B., Vassiliev invariants classify
flat braids// in S.~L.~Tabachnikov (ed.), {\em Differential and Symplectic
Topology of Knots and Curves,}
 AMS Translations Ser. 2, {\bf 190},
AMS, Providence RI, 1999, 83--102.

\bibitem{merxxx} Merkov A.B.,
Vassiliev invariants classify plane curves and doodles// Preprint, 1998, {\tt
http://www.botik.ru/\~{}duzhin/as-papers/finv-dvi.zip}

\bibitem{M} Merkov A.B., Segment--arrow diagrams and invariants
of ornaments // Matem. Sbornik {\bf 191}:11, (2000) 47--78.

\bibitem{mnev} Mnev N.E.
On manifolds of combinatorial types of projective configurations and convex
polyhedra // Soviet Math. Doklady {\bf 32} (1985), 335--337.

\bibitem{mnev2} Mnev N.E.
The universality theorems on the classification problem of configuration
varieties and convex polytope varieties, in: Topology and Geometry -- Rohlin
Seminar (O.Viro, ed.), Lecture Notes Math. {\bf 1346}, Springer-Verlag, 1988,
527--544.

\bibitem{nak} Nakamura T., The topology of the configuration of
projective subspaces in a projective space// Sci. Papers Coll. Arts Sci. Univ.
Tokyo {\bf 37} (1987), 13--35 and {\bf 41} (1991), 59--81.

\bibitem{OS} Orlik P. and Solomon L.
Combinatorics and topology of complements of hyperplanes// Invent. Math., {\bf
56}(2) (1980), 167--189.

\bibitem{OT} Orlik P. and Terao H., Arrangements of Hyperplanes,
Springer--Verlag, 1992.

\bibitem{Pham 65} Pham F., Formules de Picard--Lefschetz
g\'en\`eralis\'ees et ramification des int\'egrales// Bull. Soc. Math. France,
93 (1965), 333--367.

\bibitem{Pham 67} Pham F., Introduction \`a l'\'etude
topologique des singularit\'es de Landau, Gauthier-Villars, Paris, 1967.

\bibitem{PV} Polyak M. and Viro O. Gauss diagram formulas
for Vassiliev invariants// Internat. Math. Res. Notes {\bf 11} (1994),
445--453.

\bibitem{R} Rybnikov G.L. (1998) On the fundamental group of
the complex hyperplane arrangement// {\tt math.AG/9805056}.

\bibitem{Sal} Salvetti M. Topology of the complement of real
hyperplanes in $\C^N$// Inventiones Math. {\bf 88} (1987), 603--618.

\bibitem{Sal2} Salvetti M. The homotopy type of Artin groups.//
Math. Res. Lett., 1994, {\bf 1}:5, 565--577.

\bibitem{Schechtman-Varchenko 91} Schechtman V., Varchenko A.,
Arrangements of hyperplanes and Lie algebra homology// Invent. Math. 106
(1991), 139--194.

\bibitem{shum}
Shumakovich, A. Explicit formulas for strangeness of plane curves // Algebra i
Analiz, 1996, {\bf 7}:3, 445--472. English transl.: St. Petersburg Math. J.

\bibitem{tur1} Turchin V.E., Homology groups of complexes of
two-connected graphs //
Russian Math. Surveys, {\bf 52}:2 (1997), 426--427.

\bibitem{tur2} Turchin V.E., Homology isomorphism of the complex of
2-connected graphs and the graph complex of trees// in: Topics in Quantum
Groups and Finite-Type Invariants. Mathematics at the Independent University of
Moscow (B.~Feigin and V.~Vassiliev, eds.), AMS Translations. Ser.~2. Vol.~185.
Advances in the Mathematical Sciences. AMS, Providence RI, 1998, 145--153.

\bibitem{turdis} Tourtchine V.
Sur l'homologie des espaces des n\oe uds non-compacts, preprint 2000,
arXiv:{\tt math.QA/0010017}

\bibitem{var-} Varchenko A.N., Combinatorics and topology of dispositions
of affine hyperplanes in a real space// Funct. Anal. and its Appl., {\bf 21}:1
(1987), 11--22.

\bibitem{varchenko} Varchenko A.N., Multidimensional
hypergeometric functions in conformal field theory, algebraic K-theory,
algebraic geometry// ICM 90, p. 281--300.

\bibitem{VS} Vassiliev V.A., Serganova V.V.
On the number of real and complex moduli of singularities of smooth functions
and realizations of matroids // Math. Notes 49:1, 1991, 15--20.

\bibitem{VGZ} Vassiliev, V.A., Gelfand, I.M., and Zelevinskii, A.V.
General hypergeometric functions on complex Grassmannians // Funct. Anal. and
its Appl., {\bf 21}:1, 1987, 19--31.

\bibitem{b94} Vassiliev V.A.
Complements of discriminants of smooth maps: topology and applications, Revised
ed., Translations of Math. Monographs, V.~98, AMS, Providence RI, 1994.

\bibitem{fasis} Vassiliev V.A., Topology of complements of discriminants,
1997, 552 p. Phasis, Moscow (in Russian).

\bibitem{knsp} Vassiliev V.A. Cohomology of knot spaces//
in: Theory of Singularities and its Applications (V.~I.~Arnold, ed.), Advances
in Soviet Math. Vol. 1, AMS, Providence, RI, 1990, 23--69.

\bibitem{V-2}
Vassiliev, V.A., A geometric realization of the homology of classical Lie
groups, and complexes, $S$-dual to the flag manifolds
// St.-Petersburg Math. J., {\bf 3}:4 (1991), 809--815.

\bibitem{congr} Vassiliev V.A.
Complexes of connected graphs// in: The I.~M.~Gel'fand's mathematical seminars
1990--1992 (L.~Corvin, I.~Gel'fand, J.~Lepovsky, eds.), Birkh\"auser, Basel,
1993, 223--235.

\bibitem{ornam} Vassiliev V.A., Invariants of ornaments//
in: Singularities and Bifurcations (V.I. Arnold, ed.), Advances in Soviet
Math., {\bf 21}, AMS, Providence RI, 1994, 225--262.

\bibitem{bjo} Vassiliev V.A., Topology of two-connected graphs and
homology of spaces of knots// in: Differential and Symplectic Topology of Knots
and Curves (S.~L.~Tabachnikov, ed.), AMS Transl., Ser. 2, {\bf 190}, AMS,
Providence RI, 1999, 253--286.

\bibitem{modov} Vassiliev V.A., On finite-order invariants of
triple points free plane curves// in: D.~B.~Fuchs Anniversary Volume
(A.~Astashkevich and S.~Tabachnikov, eds.), AMS Transl., Ser. 2, {\bf 194},
AMS, Providence RI, 1999, 275--300.

\bibitem{kotor} Vassiliev V.A., Topological order complexes and
resolutions of discriminant sets// Publications de l'Institut Math\'ematique
Belgrade, Nouvelle s\'erie, {\bf 66(80)} (1999), 165--185.

\bibitem{vtt} Vassiliev V.A.,
On combinatorial formulas for cohomology of spaces of knots// Moscow Math. J.,
1:1 (2001), 91--123.

\bibitem{siersma} Vassiliev V.A.,
Resolutions of discriminants and topology of their complements// in:
Proceedings of the Advanced Summer Institute on Singularity Theory (D. Siersma,
C.T.C. Wall and V.M.Zakalyukin, eds.) Kluwer Academic Publishers, Dorderecht,
Netherlands, 2001.

\bibitem{royal} Vassiliev V.A.,
Homology of spaces of knots in any dimensions// Philosophical Trans. of the
London Royal Society. Ser. A {\bf 359}:1784(2001), 1343--1364.

\bibitem{ramif} Vassiliev, V.~A. Branching Integrals, Moscow MCCME, 2000
(in Russian).

\bibitem{Whit} Whitehead G.W., Recent Advances in Homotopy Theory.
Providence, RI: AMS, 1970, 82 p.

\bibitem{Yuz} Yuzvinsky S., Small rational model of subspace
complement// Trans. AMS, to appear. {\tt
http://xxx.lanl.gov/abs/math.CO/9806143},

\bibitem{Yuz2} Yuzvinsky S., Rational model of subspace complement on
atomic complex// Publications de l'Institut Math\'ematique Belgrade, Nouvelle
s\'erie, {\bf 66(80)} (1999), 157--164.

\bibitem{Yuz3}
Yuzvinsky S., Orlik--Solomon algebras in algebra and topology, Russian Math.
Surveys, {\bf 56}:2 (2001), 87--166.

\bibitem{Zas} Zaslavsky T., Facing up to arrangements: Face-count
formulas for partitions of space by hyperplanes// Memoirs Amer. Math. Soc.,
{\bf 1}(154), 1975.

\bibitem{Z} Ziegler G.M., On the difference between
real and complex arrangements// Math. Zeitschrift {\bf 212} (1993), 1--11.

\bibitem{Z2} Ziegler G.M., Combinatorial models for subspace arrangements,
Habilitations--Schrift, Techn. Univ. Berlin, 1992.

\bibitem{ZZ} Ziegler G.M., \v{Z}ivaljevi\'c R.T.,
Homotopy type of arrangements via diagrams of spaces// Math. Ann., {\bf 295}
(1993), 527--548.

\end{thebibliography}
\end{document}